\documentclass[reqno,11pt]{amsart}
\usepackage[margin=1in]{geometry}
\usepackage{amsfonts,amssymb,mathrsfs}
\usepackage[hidelinks]{hyperref}
\hypersetup{pdftitle={Nonlinear stability of viscous shock waves for the one-dimensional isentropic model of viscous gas with singular density-dependent viscosity},
	pdfauthor={Guanghan Zhao, Changchun Liu}}
\allowdisplaybreaks[4]
\numberwithin{equation}{section}
\newtheorem{theorem}{Theorem}[section]
\newtheorem{lemma}{Lemma}[section]
\newtheorem{proposition}{Proposition}[section]
\theoremstyle{remark}

\newcommand{\R}{\mathbb R}

\begin{document}
	\title[Nonlinear stability of viscous shock waves]
	{Nonlinear stability of viscous shock waves for the one-dimensional isentropic model of viscous gas\\ with singular density-dependent viscosity}
	\author[G. Zhao and C. Liu]
	{Guanghan Zhao and Changchun Liu$^*$}
	\address{Guanghan Zhao\hfill\break School of Mathematics, Jilin
		University, Changchun 130012, China}
	\email{zhaogh24@mails.jlu.edu.cn}
	\address{Changchun Liu\hfill\break School of Mathematics, Jilin
		University, Changchun 130012, China}
	\email{liucc@jlu.edu.cn}
	
	\thanks{This work is supported by the Science and Technology Development Plan
		Project of Jilin Province, China (No. 20250102009JC)}
	\thanks{$^*$Corresponding author.}
	\date{}
	
	\begin{abstract}
		We study the nonlinear stability of viscous shock waves for the one-dimensional isentropic Navier--Stokes equations with density-dependent viscosity $\mu(v)=av^{-m}$, where the viscous coefficient $av^{-(m+1)}$ with $a>0$ and $m>0$ possesses the strong singularity at $v=0$. Our main purpose is to establish the existence, uniqueness and nonlinear stability, up to translation, of monotone viscous shock profiles connecting a positive state to the singular end state $v_+=0$. We first determine the algebraic decay rates of the shock profiles near the singular state. Compared with the case $m=0$ in \cite{ZLLLM}, the stronger viscous singularity and slower decay of the shock profiles introduce additional difficulties in deriving the necessary energy estimates. To address the challenges, after introducing an effective velocity and using an anti-derivative method, we couple the two perturbation equation and construct a suitable weight function; thereby establish the interactive-energy-estimates to obtain the expected energy estimates for velocity. We introduce a more flexible auxiliary functional than that used in \cite{ZLLLM} to handle uncontrollable terms. Furthermore, by reformulating the perturbation system and choosing appropriate weight functions, we establish the required singular weighted energy estimates.
	\end{abstract}
	\maketitle
	
	
	\section{Introduction}
	Subsequently to the previous studies \cite{ZLLLM}, we study the Cauchy problem for the one-dimensional isentropic Navier--Stokes system with singular density-dependent diffusion in
	Lagrangian coordinates,
	\begin{equation}
		\begin{cases}
			v_t-u_x=0,\\
			u_t+p(v)_x=\left(\dfrac{\mu(v)}{v} u_x\right)_x,
		\end{cases}
		\qquad x\in\R,\quad t>0,
		\label{1-1}
	\end{equation}
	with initial data
	\begin{equation}
		(v,u)(x,0)=(v_0,u_0)(x)\longrightarrow(v_\pm,u_\pm)
		\quad\text{as }x\to\pm\infty.
		\label{1-2}
	\end{equation}
	Here $v>0$ is the specific volume, $u$ is the velocity, $p$ is the pressure, $\mu(v) = av^{-m}$ ($m > 0$) is the viscosity coefficient and $a$ is a given positive constant, which will be normalized to be $1$. It is well-known that the viscosity of the flow is not constant and depends on the temperature; see \cite{CSCG,GH}. Especially, for isentropic flows, temperature is related to density, and hence temperature-dependence viscosity provides a physical background for density-dependent viscosity laws; see \cite{KSMANT,MAWY}. In this paper, we assume that the state constants satisfy 
	\begin{equation}
		v_->v_+=0.
		\label{v-}
	\end{equation}
	Thus, the system \eqref{1-1} possesses the singularity at $v=0$. In contrast to the usual polytropic law $p(v)=av^{-\gamma}$, we consider 
	\begin{equation}
		p\in C^4([0,\infty)),\quad p'(v)<0,\quad p''(v)>0
		\quad \hbox{for all } v\geq0.
		\label{1-4}
	\end{equation}
	A typical example is $p(v)=(1+v)^{-\gamma}$ with $\gamma>0$. In particular, $p(0)$ is finite under these assumptions.
	
	\medskip
	
	{\bf Background of study.}
	When the viscosity coefficient $\mu$ is a positive constant, the nonlinear stability of viscous shock waves for system \eqref{1-1} without singular viscosity has been extensively investigated over the past several decades. Matsumura and Nishihara \cite{MANK} established the asymptotic stability of viscous shock profiles for a one-dimensional model of compressible viscous gas under sufficiently small zero-mass perturbations. Subsequently, Goodman \cite{GJ} developed the related stability theory for scalar viscous conservation laws, while Liu \cite{LT} proved the nonlinear stability of shock waves for the full compressible Navier--Stokes equations, allowing convergence to a suitably shifted shock profile. Then, Zumbrun and Howard \cite{ZKHP} developed pointwise semigroup methods and Green function estimates to link spectral stability with nonlinear stability. Within this framework, Mascia and Zumbrun \cite{MCZK} proved the nonlinear orbital stability of small-amplitude Lax shock profiles for symmetric hyperbolic-parabolic systems, including compressible Navier–Stokes equations. Moreover, under suitable spectral and structural conditions, this analysis was extended to large-amplitude shocks in \cite{MCZKS}. Spectral stability of isentropic Navier--Stokes shock waves under large-amplitude conditions was investigated in \cite{BBHJL,HJLOZK}. Recently, research has turned to the stability of composite waves. Huang and Matsumura \cite{HFMA} established the stability of a composite wave of two viscous shocks. Kang, Vasseur and Wang \cite{KMJ} proved the long-time stability of a composite wave consisting of a viscous shock and a rarefaction wave for the barotropic Navier--Stokes equations by employing the $a$-contraction method with shifts. Han, Kang and Kim \cite{HSKM} studied a composite wave composed of two small-amplitude viscous shocks, allowing their strengths to be chosen independently.  
	
	The system \eqref{1-1}, with density-dependent viscosity coefficients where the specific volume remains strictly positive, has been extensively studied in recent decades. Jiang, Xin and Zhang \cite{JXZ} established the existence of global weak solutions, while Jiu and Xin \cite{JX} further investigated the corresponding Cauchy problem. Regarding the stability of nonlinear waves, Matsumura and Wang \cite{MAWY} proved the asymptotic stability of viscous shock waves of arbitrary amplitudes for the power-law viscosity $\mu(v)=av^{-m}$ with $m\geq(\gamma-1)/2$, under sufficiently small zero-mass perturbations. Later, He, Tang and Wang \cite{HLTSWT} proved the nonlinear stability of viscous shock waves under a class of initial perturbations that allow for large density oscillations. In the same year, Vasseur and Yao \cite{VAYL} removed the restriction on the viscosity exponent and proved that, under suitable small zero-mass perturbations, viscous shock waves of arbitrary amplitude are nonlinearly stable for any $m\in\mathbb{R}$. Liu and Cui \cite{LJCM} subsequently investigated the asymptotic stability of viscous shock waves on the half-line. He and Huang \cite{HLHF} further extended the nonlinear stability theory of large-amplitude viscous shocks to general pressure laws and smooth positive viscosity coefficients. Furthermore, Huang and Yuan \cite{HFYQ} proved the stability of large-amplitude viscous shocks under small spatially periodic perturbations, while Kang and Vasseur \cite{KMJVA} used the relative entropy method to prove a contraction property for small-amplitude shocks under large perturbations. More recently, Chang, He and Ma \cite{CLHLMJ} established the stability of strong viscous shocks under small periodic perturbations in the half-space. Cui \cite{CM} studied the vanishing-viscosity limit in the presence of shock waves, proving convergence toward the corresponding Euler shock. 
	
	Recently, considerable progress has been made in the stability analysis of nonlinear waves for Burgers equations with singular diffusion of the form
	\begin{equation*}
		u_t+f(u)_x=\mu\left(\frac{u_x}{u^{m+1}}\right)_x,\qquad m>-1.
	\end{equation*}
	The existence and nonlinear stability of viscous shock profiles have been investigated in the critical fast-diffusion case $m=0$ by Li, Li, Mei and Nave \cite{Li-Li-Mei-Nave}; in the fast-diffusion case $-1<m<0$, including both nondegenerate and degenerate shocks, by Xu, Mei, Nave and Sheng \cite{Xu-Mei-Nave-Sheng}; in the super-fast-diffusion case $m>0$ by Li, Li and Mei \cite{Li-Li-Mei}.
	
	Furthermore, Zhao et al. \cite{ZLLLM} established the existence, uniqueness and nonlinear asymptotic stability of viscous shock waves connecting a positive state to the singular state $v_+=0$ for the system \eqref{1-1} in the case $m=0$.
	
	\medskip
	
	{\bf Main difficulty and strategy.}
	In this paper, we further extend the study to the case $m>0$. We aim to achieve two main objectives. First, we establish the existence and uniqueness, up to translation, of a monotone viscous shock profile $(V,U)(x-st)$ ($s$ is the shock speed) connecting $(v_-,u_-)$ and $(0,u_+)$, and show its exact decay rates as $x-st \to \pm \infty$. Second, we prove the nonlinear asymptotic stability of the viscous shock profile $(V,U)$.
	
	For the nonlinear stability of $(V,U)$, in the classical setting where $v_{+},v_{-}>0$, a common approach is to introduce the effective velocity 
	\begin{equation*}
		b=u-\mu v_x/v,
	\end{equation*}
	and first establish the nonlinear stability of the transformed shock profile $(V,B)$ before recovering the stability of $(V,U)$. However, due to the presence of singularity, we directly investigate the nonlinear stability of the original shock profile $(V,U)$ after introducing the effective velocity $b$. 
	
	Compared with the case $m=0$ studied in \cite{ZLLLM}, the main difficulty is that the viscosity coefficient $v^{-(m+1)}$ exhibits a stronger singularity as $v \to 0^+$ for $m>0$. Moreover, the shock profile shows a slower algebraic decay toward the singular end state. This makes the nonlinear terms of the perturbation equations more difficult to control. Consequently, the method developed for $m=0$ cannot be applied directly. 
	
	To derive a weighted space-time $L^{2}$ bound of the velocity perturbation $q$, we exploit the coupling between the two perturbation equations to establish interactive-energy-estimates, where a new weight function adapted to asymptotic behavior of the shock profile is required. Moreover, to control certain terms that cannot be directly handled by the energy estimates, we introduce a more flexible auxiliary functional and establish the corresponding estimates. We then exploit the relation between the effective velocity perturbation $\psi_z$ and $q$ to reformulate the perturbation system solely in terms of the specific-volume perturbation $\phi_z$ and $q$. Similarly, based on this reformulation, we construct appropriate weight functions, and then derive the singular weighted-energy estimates to establish the nonlinear stability of $(V,U)$. Finally, a cut-off technique is employed to handle the boundary terms arising at spatial infinity and complete the proof.

	\medskip
	Before proceeding, we introduce some notations for convenience.
	
	\noindent\textbf{Notation.}
	Let $C$ be a generic positive constant that may vary between lines. The integrals $\int_{\mathbb{R}} f(x) \mathrm{d} x$ and $\int_0^t \int_{\mathbb{R}} f(x, \tau) \mathrm{d} x \mathrm{d}r$ will be expressed as $\int f(x)$ and $\int_0^t \int f(x, \tau)$, respectively.
	We denote $\|f\|:=\|f\|_{L^2(\R)}$ and
	$\|f\|_k:=\|f\|_{H^k(\R)}$ $(k > 0)$. For a positive weight $w$, set
	\[
	\|f\|_w^2=\int w|f|^2,\qquad
	\|f\|_{k,w}^2=\sum_{j=0}^k\int w|\partial_z^j f|^2.
	\]
	The corresponding spaces are denoted by $L_w^2$ and $H_w^k$.
	For a vector, $\| (f,g,h) \|^{2}$ is equivalent to $\| f \|^{2} + \| g \|^{2} + \| h \|^{2}$.
	The notation $f\sim g$ means $f/g\to1$ at the specified limit.
	
	An admissible viscous shock profile is a nonconstant pair
	$(V,U)\in C^4(\R)$, with $V>0$, such that
	\begin{equation}
		(v,u)(x,t)=(V,U)(z),\qquad z=x-st,
		\label{1-111}
	\end{equation}
	solves \eqref{1-1}. Substituting it into \eqref{1-1}, we get the equation of shock profile
	\begin{equation}
		\begin{cases}
			-sV_z=U_z,\\
			-sU_z+p(V)_z=\left(\dfrac{U_z}{V^{m+1}}\right)_z,
		\end{cases}
		\label{1-5}
	\end{equation}
	which satisfies
	\begin{equation}
		\begin{gathered}
			V(-\infty)=v_-,\qquad V(+\infty)=0,\qquad U(\pm\infty)=u_\pm,\\
			V_z(\pm\infty)=U_z(\pm\infty)=0,\qquad
			\lim_{z\to\pm\infty}\frac{U_z(z)}{(V(z))^{m+1}}=0.
		\end{gathered}
		\label{1-6}
	\end{equation}
	The last condition is a vanishing viscous-stress condition. 
	
	Under the conditions of \eqref{1-6}, we first state the sufficient and necessary conditions for the existence of viscous shock profile $(V,U)$ to system \eqref{1-1} as follows.
	
	\begin{theorem}[Existence and uniqueness of the profile]\label{theorem1.1}
		Assume that \eqref{1-4} holds. The Rankine--Hugoniot conditions for an admissible viscous shock profile connecting $(v_{-},u_{-})$ to $(0,u_{+})$ are
		\begin{equation}
			u_+-u_-=sv_-,\qquad s^2=\frac{p(0)-p(v_-)}{v_-}.
			\label{1-7}
		\end{equation}
		The connecting profile necessarily has $s<0$. Strict convexity of $p$ gives the
		Lax's entropy condition
		\begin{equation}
			-\sqrt{-p'(0)}<s<-\sqrt{-p'(v_-)}.
			\label{1-8}
		\end{equation}
		The viscous shock profile satisfying \eqref{1-6} exists
		if and only if \eqref{1-7} and \eqref{1-8} hold with $s<0$. The profile is unique up to a shift and satisfies
		\begin{equation}
			V_z<0,\qquad U_z<0\qquad\text{on }\R.
			\label{1-9}
		\end{equation}
		Moreover, 
		\begin{equation}
			V(z) \sim (a_+z)^{-\frac{1}{m+1}} \quad(z\to+\infty),\qquad
			v_--V(z)\sim A_-e^{\lambda_-z}\quad(z\to-\infty),
			\label{1-911}
		\end{equation}
		where
		\begin{equation*}
			a_+=\frac{(p'(0)+s^2)(m+1)}{s}>0, \qquad \lambda_-=-\frac{v_-^{m+1}}{s} (p'(v_-)+s^2)>0,
		\end{equation*}
		and $A_->0$ is a constant depending on the translation.
	\end{theorem}
	
	In order to deal with the nonlinear singular viscous term $(\frac{u_{x}}{v^{m+1}})_{x}$ better, we introduce the effective velocity
	\begin{equation}
		b = u - \frac{v_x}{v^{m+1}}.
		\label{1-10}
	\end{equation}
	Then, $(v,b)$ satisfies
	\begin{equation}
		\begin{cases}
			v_t-b_x=\left(\dfrac{v_x}{v^{m+1}}\right)_x,\\
			b_t+p(v)_x=0,
		\end{cases}
		\label{1-11}
	\end{equation}
	with the initial data
	\begin{equation}
		(v,b)(x,0) = (v_0,b_0)(x) \to (v_\pm,b_\pm), \quad \text{as } x \to \pm \infty,
		\label{1-12}
	\end{equation}
	where
	\begin{equation*}
		b_0 = u_0 - \frac{v_{0x}}{v_0^{m+1}}. 
	\end{equation*}
	Furthermore, the stress condition in \eqref{1-6} gives $b_\pm=u_\pm$.
	
	Effective velocity formulations for related compressible models have been applied in \cite{HLHF, KMJVA, VAYL}; for viscoelastic systems, see also \cite{LJGZ}. The relevant structures of weak solutions also appear in \cite{BDBD,BDBDO}.
	
	\begin{theorem}\label{theorem1.2}
		Let $(V,U)$ be the profile in Theorem \ref{theorem1.1}, and define
		\[
		B = U- \frac{V_z}{V^{m+1}}.
		\]
		Then, the system \eqref{1-11} has a unique (up to a shift) smooth monotone traveling wave solution $(V,B)(x-st)$, with the end states
		$(v_-,u_-)$ and $(0,u_+)$, satisfying
		\begin{equation*}
			V_z<0,\qquad B_z<0 \quad\text{on }\;\R.
		\end{equation*}
	\end{theorem}
	Assume that the initial perturbations satisfy 
	\begin{equation*}
		v_0-V,\quad b_0-B\in L^1(\mathbb R).
	\end{equation*}
	We choose a shift $x_0 \in \mathbb R$ satisfying the compatibility condition 
	\begin{equation}
		\int\binom{v_0(x)-V(x)}{b_0(x)-B(x)} \mathrm{d}x = x_0\binom{-v_-}{b_+-b_-}.
		\label{3-1}
	\end{equation}
	For this choice of $x_0$, we define 
	\begin{equation}
		V_0(x):=V(x+x_0), \quad U_0(x):=U(x+x_0), \quad B_0(x):=B(x+x_0). 
		\label{1-151}
	\end{equation}
	Using the end states of $V$ and $B$, we have
	\begin{equation*}
		\int \binom{V(x)-V_{0}(x)}{B(x)-B_{0}(x)} \mathrm{d} x = - x_0\binom{-v_{-}}{b_{+}-b_{-}}.
	\end{equation*}
	Consequently, combining this identity with \eqref{3-1}, we obtain the zero-mass conditions
	\begin{equation}
		\int \binom{v_{0}(x)-V_{0}(x)}{b_{0}(x)-B_{0}(x)} \mathrm{d} x = \int \binom{v_0(x)-V(x)}{b_0(x)-B(x)} \mathrm{d} x + \int \binom{V(x)-V_{0}(x)}{B(x)-B_{0}(x)} \mathrm{d} x = \binom{0}{0}.
		\label{3-2}
	\end{equation}
	Let us define
	\begin{equation*}
		(\phi_0,\psi_0)(x)=\int_{-\infty}^x
		(v_0(y)-V_0(y),\ b_0(y)-B_0(y))\,\mathrm{d}y.
	\end{equation*}
	Then $(\phi_0,\psi_0)(\pm\infty)=0$. We are now state the nonlinear stability results for the viscous shock profile $(V,U)(x-st)$ of \eqref{1-1}--\eqref{1-2}.
	
	\begin{theorem}\label{theorem1.3}
		Let $m>0$ and $(V,B)(x-st+x_{0})$ be the shock profile obtained in Theorem \ref{theorem1.2}. Assume that $v_{0} - V_{0} \in L^{1}(\mathbb R)$, $b_{0} - B_{0} \in L^{1}(\mathbb R)$ and \eqref{3-1} holds. $V_{0}$, $U_{0}$ and $B_{0}$ are defined in \eqref{1-151}. Set $w_{2,0}=V_0^{-2(m+1)}$. There exists $\delta>0$
		such that, if
		\begin{align*}
			\|\phi_0\| + \|\psi_0\|_2 + \|v_0-V_0\|_{1,w_{2,0}} + \|u_0-U_0\|_{1,w_{2,0}}\leq\delta,
		\end{align*}
		then \eqref{1-1}--\eqref{1-2} has a unique global solution $(v,u)(x,t)$ satisfying
		\[
		v-V \in C([0,\infty);H^1_{w_{2}}),\quad u-U \in C([0,\infty);H^1_{w_{2}}),
		\]
		where $w_{2} = V^{-2(m+1)}$. Moreover, the solution has the following asymptotic convergence
		\begin{equation*}
			\sup_{x\in\R}|(v,u)(x,t)-(V,U)(x-st+x_0)|\to 0,
			\quad\text{as }\;t\to\infty.
		\end{equation*}
	\end{theorem}
	The rest of this paper is organized as follows. In Section~\ref{sec:profiles}, the existence, uniqueness, monotonicity and asymptotic behavior of the shock profile $(V,U)$ will be discussed, and then of $(V,B)$. In Section~\ref{sec:stability}, we show the details of artful weighted energy estimates and derive the nonlinear stability results.

	\section{Existence of shock profiles}\label{sec:profiles}
	\begin{proof}[Proof of Theorem \ref{theorem1.1}]
		Suppose first that an admissible profile $(V,U)(x-st)$ exists. Integrating the first
		equation of \eqref{1-5} gives
		\begin{equation*}
			u_+-u_-=sv_-.
		\end{equation*}
		Substituting \eqref{1-5}$_{1}$ into \eqref{1-5}$_{2}$ yields
		\begin{equation}
			s^2V_z+p(V)_z=-s \left(\frac{V_z}{V^{m+1}}\right)_z.
			\label{2-1}
		\end{equation}
		By $-sV_z=U_z$, the last condition of \eqref{1-6} leads to $V_z/V^{m+1}\to0$ as $z \to \pm \infty$.
		We integrate \eqref{2-1} over $\mathbb{R}$ to get
		\begin{equation}
			s^2=\frac{p(0)-p(v_-)}{v_-}.
			\label{2-2}
		\end{equation}
		Integrating \eqref{2-1} from $z$ to $+\infty$, we obtain
		\begin{equation*}
			s^2V+p(V)+s \frac{V_z}{V^{m+1}}=p(0). 
		\end{equation*}
		Thus, it follows that
		\begin{equation}
			f(V):=V_z=-\frac{V^{m+1}}{s}g(V),\qquad
			g(V):=p(V)-p(0)+s^2V.
			\label{2-4}
		\end{equation}
		Since $g''(V)=p''(V)>0$ and $g(0)=g(v_-)=0$, by strict convexity,
		\begin{equation*}
			g(V)<0, \quad \hbox{for } 0<V<v_-. 
		\end{equation*}
		Uniqueness for the scalar autonomous equation prevents a nonconstant connecting orbit from crossing either equilibrium. The prescribed end states then imply $V_{z} < 0$. Since $V>0$, by \eqref{2-4}, we get
		\begin{equation}
			s<0.
			\label{2-6}
		\end{equation}
		Using strict convexity of $p$, we get
		\begin{equation*}
			p'(0)<\frac{p(v_-)-p(0)}{v_-}<p'(v_-).
		\end{equation*}
		In view of \eqref{2-2} and \eqref{2-6}, we have \eqref{1-8}.
		
		Conversely, suppose that \eqref{1-7} and \eqref{1-8} hold with $s<0$.
		The definition of $f$ in \eqref{2-4}, $s<0$ and $g(V)<0$ lead to $f<0$ on $(0,v_-)$. We only need to find the global solution of the following ODE
		\begin{align}
			\frac{\mathrm{d}z}{\mathrm{d}V} = \frac{1}{f(V)}.
			\label{2-81}
		\end{align}
		Hence, for any $v_*\in(0,v_-)$ and $z_*\in\R$ satisfying $V(z_*)=v_*$ and $z(v_*)=z_*$, integrating \eqref{2-81} yields
		\begin{equation*}
			z-z_*=F(V):=\int_{v_*}^V\frac{\mathrm{d}y}{f(y)}.
		\end{equation*}
		As $F'(V)=1/f(V)<0$, the function $F$ is finite and monotonically decreasing for any given $V \in (0,v_-)$. By Taylor's expansion of \eqref{2-4} at $V=0$, we derive
		\begin{equation}
			f(V)=- \frac{V^{m+2}}{s} \left[\frac{p(V)-p(0)}{V} + s^2\right] = -\frac{a_+}{m+1}V^{m+2}+O(V^{m+3}).
			\label{2-10}
		\end{equation}
		From \eqref{2-2}, at $V=v_{-}$, we can also get 
		\begin{equation}
			\begin{aligned}
				f(V) = - \frac{V^{m+1}}{s} [p(V) - p(v_{-}) + s^2 (V - v_{-})] = & - \frac{V^{m+1} (V - v_{-})}{s} \left[\frac{p(V)-p(v_{-})}{V - v_{-}} + s^2\right] 
				\\
				= & \lambda_-(V-v_-)+O((V-v_-)^2).
			\end{aligned}
			\label{2-11}
		\end{equation}
		Consequently, it becomes apparent that
		\begin{equation*}
			\lim_{V \to 0^+} F(V) = +\infty,
			\quad
			\lim_{V \to v_-} F(V) = -\infty. 
		\end{equation*}
		Then $V$ is a strictly decreasing solution of \eqref{2-4} with the limiting values $V(\pm \infty) = v_{\pm}$. Therefore, according to the inverse function theorem, there exists a unique continuous inverse function $F^{-1}$ such that
		\begin{equation*}
			V = F^{-1}(z-z_\star), \quad z\in(-\infty,+\infty).
		\end{equation*}
		Integrating \eqref{1-5}$_{1}$ over $(z,+\infty)$, the shock wave profile $U(z)$ is uniquely determined. Furthermore, by virtue of $U_z = - s V_z < 0$, \eqref{1-9} holds.
		In particular, \eqref{2-4} and $g(0) = g(v_-) = 0$ imply
		\begin{align*}
			\frac{V_z}{V^{m+1}} = - \frac{g(V)}{s} \to 0,
			\quad \text{as } z \to \pm \infty.
		\end{align*}
		Thus, all conditions in \eqref{1-6} hold. 
		
		Next, we are going to show the convergence rates for $V(z)$ as $z \to \pm \infty$. Using \eqref{2-10}--\eqref{2-11} and putting $Y=v_- - V$ gives
		\begin{equation*}
			\begin{cases}
				-\dfrac{(m+1)V_{z}}{V^{m+2}}=\left(\dfrac{1}{V^{m+1}}\right)_z=a_++O(V),
				\\
				Y_z=\lambda_-Y+O(Y^2),
			\end{cases}
		\end{equation*}
		which implies $V(z) \sim (a_{+}z)^{-\frac{1}{m+1}} \to 0$ as $z \to +\infty$ and $Y$ decays exponentially as $z\to-\infty$ by comparison. Note that $Y_z/Y-\lambda_-=(\log Y - \lambda_- z)_z=O(Y)$ is integrable. Integrating it over $[z,z_{0}]$ ($z_{0}\in\mathbb{R}^{+}$) and then exponentiating yields
		\begin{equation*}
			Y(z)e^{-\lambda_-z}=Y(z_{0})e^{-\lambda_-z_{0}}\exp\left\{- \int_{z}^{z_{0}} O(Y)\right\}.
		\end{equation*}
		Thus, $Y(z)e^{-\lambda_-z}$ tends to a constant $A_- > 0$ as $z \to -\infty$. The proof of \eqref{1-911} is complete.
	\end{proof}
	
	\begin{proof}[Proof of Theorem \ref{theorem1.2}]
		Define
		\begin{equation}
			B=U-\mu\frac{V_z}{V}.
			\label{2-15}
		\end{equation}
		Using \eqref{1-5}, we obtain
		\begin{equation}
			\begin{cases}
				-sV_z-B_z=\mu\left(\dfrac{V_z}{V}\right)_z,\\
				-sB_z+p(V)_z=0.
			\end{cases}
			\label{2-17}
		\end{equation}
		Since $p'(V)<0$, $s<0$ and $V_z<0$, we conclude that $B_z<0$. Uniqueness up to translation follows from the uniqueness of $V$ and the definition \eqref{2-15}.
	\end{proof}
	
	We first establish some basic properties of the shock profile $V(z)$.
	\begin{lemma}\label{lemma3.1}
		Let $V(z)$ be the shock wave obtained in Theorem \ref{theorem1.1}. It holds that
		\begin{equation}
			|V_z| \leq C V^{m+2}, \quad |V_{z z}| \leq C V^{2m+3}, \quad |V_{zzz}| \leq C V^{3m+4}.
			\label{3-16}
		\end{equation}
		Moreover,
		\begin{equation}
			\left|\frac{V_z}{V^{m+2}}\right|\leq C,\quad
			\left|\left(\frac{V_z}{V^{m+2}}\right)_z\right|
			+\left|\left(\frac{V_z}{-p'(V)V^{m+2}}\right)_z\right|
			\leq C|V_z|,
			\label{profile:coeff}
		\end{equation}
		for any $z \in(-\infty,+\infty)$.
	\end{lemma}
	\begin{proof}
		From \eqref{2-4} and the fact that $0<V \leq v_-$, there exists a constant $C > 0$ such that $|V_z|\leq CV^{m+2}$. Moreover, since $|f'(V)|\leq CV^{m+1}$ and $|f''(V)|\leq CV^{m}$, differentiating the profile equation yields
		\begin{equation*}
			|V_{zz}| = |f'(V)V_{z}| \leq CV^{m+1}|V_z| \leq C V^{2m+3},
		\end{equation*}
		and
		\begin{equation*}
			|V_{zzz}| = |f'(V)V_{zz}+f''(V)V_{z}^{2}| \leq C(V^{m+1} |V_{zz}| + V^{m} V_z^{2}) \leq C V^{3m+4}.
		\end{equation*}
		Next, we establish estimates of \eqref{profile:coeff}. In view of \eqref{2-4}, we have
		\begin{equation*}
			\frac{V_z}{V^{m+2}} = - \frac{1}{s} \left(\frac{p(V) - p(0)}{V} + s^2\right) =-\frac{1}{s}\left(\int_0^1p'(\theta V) \mathrm{d}\theta+s^2\right). 
		\end{equation*}
		Consequently, $|V_{z}/V^{m+2}| \leq C$. Furthermore, using the regularity assumption $p \in C^{4}([0,\infty))$ and differentiating the above identity with respect to $z$, it holds that
		\begin{equation}
			\begin{aligned}
				\left|\left(\frac{V_{z}}{V^{m+2}}\right)_{z}\right| = \left|\frac{V_{z}}{s} \frac{\mathrm{d}}{\mathrm{d}V} \left(\int_0^1 p'(\theta V) \mathrm{d}\theta+s^2\right)\right| = \left|\frac{V_{z}}{s} \int_0^1 \theta p''(\theta V) \mathrm{d}\theta\right| \leq C |V_{z}|.
			\end{aligned}
			\label{3.23}
		\end{equation}
		Finally, combining \eqref{3-16} with \eqref{3.23}, we deduce that
		\begin{equation}
			\left|\left(\frac{V_{z}}{-p'(V) V^{m+2}}\right)_{z}\right| = \left|\frac{1}{-p'(V)} \left(\frac{V_{z}}{V^{m+2}}\right)_{z} + \frac{V_{z}}{V^{m+2}} \frac{p''(V) V_{z}}{(p'(V))^{2}}\right| \leq C |V_{z}|.
			\label{3.24}
		\end{equation}
		This completes the proof.
	\end{proof}

	\section{Nonlinear asymptotic stability}\label{sec:stability}
	In this section, we prove Theorem \ref{theorem1.3}, thereby establish the nonlinear stability of the shock profile for the system \eqref{1-1}--\eqref{1-2}. Since the viscosity coefficient $\mu(v)$ becomes singular as $z\to+\infty$, we shall devise new strategies to address this challenge. As a preliminary step, we use the anti-derivative technique and reformulate the problem.
	
	\subsection{Reformulation of the problem and the main results}
	Let $(V,B)(x-st)$ be the viscous shock profile obtained in Theorem \ref{theorem1.2}. By the conservation of mass and \eqref{3-2}, we have
	\begin{equation}
		\int \binom{v(x,t)-V(x+x_{0}-st)}{b(x,t)-B(x+x_{0}-st)} \mathrm{d}x
		= \int \binom{v_{0}(x)-V_{0}(x)}{b_{0}(x)-B_{0}(x)} \mathrm{d}x
		= \binom{0}{0}.
		\label{3-21}
	\end{equation}
	We thus decompose the solution of \eqref{1-1} as
	\begin{align}
		(v, b)(x, t)=(V, B)\left(x+x_0-s t\right)+\left(\phi_z, \psi_z\right)(z, t),
		\label{3-5}
	\end{align}
	where $z = x - s t$, that is, for all $z \in \mathbb{R}$ and $t \geq 0$, 
	\begin{equation*}
		(\phi, \psi)(z, t)=\int_{-\infty}^z\left(v(y + st, t)-V\left(y+x_0\right), b(y + st, t)-B\left(y+x_0\right)\right) \mathrm{d}y. 
	\end{equation*}
	By \eqref{3-21}, we obtain
	\begin{equation}
		(\phi, \psi)( \pm \infty, t)=0, \quad \hbox{for all } t \geq 0.
		\label{3-4}
	\end{equation}
	Substituting \eqref{3-5} into \eqref{1-11}, making use of \eqref{2-17} and \eqref{3-4}, and integrating the resulting equations with respect to $z$, the problem \eqref{1-1}--\eqref{1-2} is reduced to 
	\begin{equation}
		\begin{cases}
			\phi_t = s \phi_z + \psi_z + \left(\frac{V_z+\phi_{z z}}{(V+\phi_z)^{m+1}} - \frac{V_z}{V^{m+1}}\right),
			\\
			\psi_t = s \psi_z - (p(V+\phi_z)-p(V)).
		\end{cases}
		\label{3-611}
	\end{equation}
	Without loss of generality, we normalize the translation by taking $x_0=0$. The initial data are
	\begin{align}
		(\phi_0, \psi_0)(z) := (\phi, \psi)(z, 0)=\int_{-\infty}^z (v_0-V, b_0-B)(y) \mathrm{d}y,
		\label{3-9}
	\end{align}
	with $(\phi_0, \psi_0)( \pm \infty)=0$ due to \eqref{3-2}.
	
	For convenience, we introduce the perturbation variables
	\begin{equation}
		r:=\phi_{z}, \quad q:=\phi_t-s\phi_z=\psi_z + \left(\frac{V_z+\phi_{z z}}{(V+\phi_z)^{m+1}} - \frac{V_z}{V^{m+1}}\right).
		\label{defs:r-q}
	\end{equation}
	It follows from \eqref{1-10} and \eqref{2-15} that $q$ represents the perturbation of the original velocity $u$. Accordingly, the initial values of $r$ and $q$ are given by
	\begin{equation}
		r_{0}=\phi_{0z}, \quad q_0=\psi_{0z} + \left(\frac{V_z+\phi_{0z z}}{(V+\phi_{0z})^{m+1}} - \frac{V_z}{V^{m+1}}\right).
		\label{q:init}
	\end{equation}
	To handle the singular behavior of the shock profile near the zero end state, we introduce the weight functions
	\begin{equation*}
		w_i=V^{-i(m+1)},\quad i=1,2,3. 
	\end{equation*}
	For any fixed $T>0$, we define the solution space
	\begin{equation*}
		X(0,T)=\{(\phi,\psi):\ \phi\in C([0,T];L^2),\ \psi\in C([0,T];H^2),\ r,q\in C([0,T];H^1_{w_2})\}.
	\end{equation*}
	Furthermore, we introduce the energy functional 
	\begin{equation}
		\mathcal E(t)=\|\phi(t)\|^2+\|\psi(t)\|_2^2
		+\|r(t)\|_{1,w_2}^2+\|q(t)\|_{1,w_2}^2,\quad
		N(t)=\sup_{0\leq\tau\leq t}\mathcal E(\tau)^{1/2}.
		\label{3-11}
	\end{equation}
	The local existence result can be established by following the approach in \cite{ZLLLM}.
	\begin{proposition}[Local existence]\label{proposition3.1}
		Suppose that $\phi_0\in L^2(\mathbb{R})$, $\psi_0\in H^2(\mathbb{R})$ and
		$r_{0},q_0\in H^1_{w_2}(\mathbb{R})$, with $r_{0}$ and $q_0$ given by \eqref{q:init}. For any sufficiently small $\delta_0 > 0$, there exists a positive constant $T_0$ depending on $\delta_0$ such that if $N(0) \leq \delta_0$, then the problem \eqref{3-611} with \eqref{3-9} has a unique solution $(\phi, \psi) \in X(0, T_0)$ satisfying $N(t) \leq C_{*} N(0)$, for any $0 \leq t \leq T_0$, where $C_{*}$ is independent of $t$.
	\end{proposition}
	Before deriving the a priori estimates, we establish several preliminary bounds under the smallness assumption of $N(t)$. The one-dimensional Sobolev inequality $\|f\|_{L^{\infty}} \leq \sqrt{2} \|f\|^{1/2}\|f_{z}\|^{1/2}$ and \eqref{3-11} yield
	\begin{equation}
		\sup _{\tau \in[0, t]}\left\{\|(\phi, \psi, \psi_z)(\cdot, \tau)\|_{L^{\infty}}, \left\|(r/V^{m+1}, q/V^{m+1})(\cdot, \tau)\right\|_{L^{\infty}}\right\} \leq C N(t),
		\label{3-15}
	\end{equation}
	where $C>0$ is independent of $t$ and $T$. Since $0 < V < v_-$ and $m>0$, it follows that
	\begin{equation}
		|r/V| \leq C|r/V^{m+1}| \leq CN(t) \ll 1, \quad 1 + r/V \geq 1 - CN(t) > 0, \quad V + r \geq V - CN(t) V > 0,
		\label{3-151}
	\end{equation}
	provided that $N(t)$ is sufficiently small. On the other hand, using the definition of $q$ in \eqref{defs:r-q}, we have
	\begin{equation}
		\begin{aligned}
			\frac{r_{z}}{V^{m+1}} = & \left(\frac{V+r}{V}\right)^{m+1} (q-\psi_z) + \frac{V_z(V+r)^{m+1}}{V^{2(m+1)}} - \frac{V_z}{V^{m+1}} 
			\\
			= & \left(1+ \frac{r}{V}\right)^{m+1} (q-\psi_z) + \frac{V_z [(V+r)^{m+1} - V^{m+1}]}{V^{2(m+1)}}.
		\end{aligned}
		\label{rz:identity}
	\end{equation}
	By \eqref{3-151} and the mean value theorem, we have $|(V+r)^{m+1} - V^{m+1}| \leq C V^{m}|r|$. Combining this inequality with \eqref{3-16} and \eqref{3-15}, we deduce that 
	\begin{equation}
		\|r_{z}/V^{m+1}\|_{L^\infty}\leq CN(t).
		\label{rz:infinity}
	\end{equation}
	To handle the nonlinear viscous term, we introduce the auxiliary variable
	\begin{equation}
		\eta = -m^{-1}[(V+r)^{-m} - V^{-m}].
		\label{3-6111}
	\end{equation}
	Differentiating it with respect to $z$, we find
	\begin{equation}
		\eta_{z} = \frac{V_z+r_{z}}{(V+r)^{m+1}} - \frac{V_z}{V^{m+1}}. 
		\label{3-6112}
	\end{equation}
	Consequently, the velocity perturbation can be expressed as
	\begin{equation}
		q=\phi_t-s\phi_z=\psi_z+\eta_z.
		\label{3-61}
	\end{equation}
	Using this identity, we reformulate the system \eqref{3-611} as
	\begin{equation}
		\begin{cases}
			\phi_t=sr+\psi_z+(r/V^{m+1})_z+W_z,
			\\
			\psi_t=s\psi_z-p'(V)r-G,
		\end{cases}
		\label{3-6}
	\end{equation}
	where
	\begin{equation}
		W := -m^{-1}[(V+r)^{-m} - V^{-m} + m r/V^{m+1}],
		\label{3-7}
	\end{equation}
	and
	\begin{equation}
		G := p(V+r)-p(V)-p'(V)r.
		\label{3-8}
	\end{equation}
	We next define
	\begin{equation*}
		d(V) := \frac{s p''(V)V_z}{p'(V)^2}.
	\end{equation*}
	Since $p \in C^{4}([0,\infty))$, we get $d(V) \sim C|V_z|$. We denote the weighted energy estimate used below as 
	\begin{equation*}
		\mathcal D(t)=\|r(t)\|_{2,w_1}^2+\|q(t)\|_{w_1}^2
		+\|q_z(t)\|_{1,w_3}^2+
		\int d(V)(\psi^2+\psi_z^2+\psi_{zz}^2)(t).
	\end{equation*}
	With these preparations, we state the a priori estimate.
	\begin{proposition}[A priori estimate]\label{proposition3.2}
		Let the assumptions of Proposition \ref{proposition3.1} hold. Assume that $(\phi, \psi) \in X(0, T)$ is a solution of \eqref{3-611} with \eqref{3-9} for a positive constant $T$. Then there exist a constant $\delta_{1}$ satisfying $0 < \delta_{1} < \delta_{0}$, independent of $T$, such that if
		\begin{equation}
			N(t) \leq \delta_1,
			\label{3-14}
		\end{equation}
		for any $0 \leq t \leq T$, it holds that
		\begin{equation}
			\mathcal E(t)+\int_0^t\mathcal D(\tau) \mathrm{d}\tau
			\leq C\mathcal E(0).
			\label{3-12}
		\end{equation}
	\end{proposition}
	Combining the local existence result with the a priori estimate and a standard continuation argument, we obtain the following global stability theorem.
	\begin{theorem}[Global stability of the perturbation problem]\label{theorem3.1}
		Under the assumptions of Proposition \ref{proposition3.1}, there exists
		a constant $\delta_2>0$ such that if $N(0)\leq\delta_2$, the Cauchy problem \eqref{3-611} with \eqref{3-9} has a unique global solution $(\phi, \psi) \in X(0, \infty)$ satisfying \eqref{3-12} for any $t \in [0, \infty)$. Moreover, it follows that
		\begin{equation}
			\|r(t)\|_{L^\infty}+\|q(t)\|_{L^\infty} \to 0,
			\quad\text{as }t\to\infty.
			\label{3-13}
		\end{equation}
	\end{theorem}

	\subsection{Useful inequalities, basic and interactive energy estimates}
	\begin{lemma}\label{lemma3.2}
		Under the a priori assumption \eqref{3-14}, if $N(T)$ is sufficiently small, then it holds that
		\begin{equation}
			\begin{aligned}
				|W| & \leq C r^{2}/V^{m+2}, \quad \quad \quad \quad \quad \quad \quad \quad \; |G| \leq C r^2,\\
				|W_z| & \leq C \left(\frac{|V_z|}{V^{m+3}} r^2+\frac{|r||r_{z}|}{V^{m+2}}\right), \quad \quad \quad |G_z| \leq C (|V_z| r^2+|r||r_{z}|),\\
				|W_{zz}|&\leq C\left[\left(\frac{V_z^2}{V^{m+4}} + \frac{|V_{zz}|}{V^{m+3}}\right)r^2 + \frac{r_{z}^2}{V^{m+2}} + \frac{|V_z|}{V^{m+3}}|r||r_{z}| + \frac{|r||r_{zz}|}{V^{m+2}}\right],\\
				|G_{zz}|&\leq C (V_z^2r^2 + |V_z||r||r_{z}| + |V_{zz}|r^2 + r_{z}^2 + |r||r_{zz}|).
			\end{aligned}
			\label{3-18}
		\end{equation}
	\end{lemma}
	\begin{proof}
		By the definitions of $W$ and $G$ in \eqref{3-7} and \eqref{3-8}, a direct calculation gives $|W| \leq C r^{2}/V^{m+2}$ and $|G| \leq C r^2$. Furthermore, it holds that
		\begin{align*}
			W_z = & \frac{V_z+r_{z}}{(V+r)^{m+1}} - \frac{V_z}{V^{m+1}} -\frac{r_{z}}{V^{m+1}}+\frac{(m+1)V_z r}{V^{m+2}} 
			\\
			= & \left(\frac{1}{(V+r)^{m+1}}-\frac{1}{V^{m+1}}+\frac{(m+1)r}{V^{m+2}}\right) V_z + \left(\frac{1}{(V+r)^{m+1}}-\frac{1}{V^{m+1}}\right) r_{z},
		\end{align*}
		and
		\begin{align*}
			G_z & =p'(V+r)(V_z+r_{z})-p'(V) V_z-p''(V) V_z r-p'(V) r_{z}
			\\
			& =(p'(V+r)-p'(V)-p''(V) r) V_z+(p'(V+r)-p'(V)) r_{z}.
		\end{align*}
		Taylor's expansion then yields the estimates of $|W_z|$ and $|G_z|$. An analogous argument gives
		\begin{align*}
			W_{zz} = & \left(\frac{V_z+r_{z}}{(V+r)^{m+1}} - \frac{V_z}{V^{m+1}} -\frac{r_{z}}{V^{m+1}}+\frac{(m+1)V_z r}{V^{m+2}}\right)_z
			\\
			= & \left[\frac{1}{(V+r)^{m+1}}-\frac{1}{V^{m+1}}+\frac{(m+1)r}{V^{m+2}}\right]V_{zz} - (m+1) \left[\frac{1}{(V+r)^{m+2}} - \frac{1}{V^{m+2}} + \frac{(m+2)r}{V^{m+3}}\right]V_z^2
			\\
			& - \frac{(m+1)r_{z}^2}{(V+r)^{m+2}} - 2 (m+1)\left(\frac{1}{(V+r)^{m+2}} -  \frac{1}{V^{m+2}}\right)V_z r_{z} + \left(\frac{1}{(V+r)^{m+1}} - \frac{1}{V^{m+1}}\right)r_{zz},
		\end{align*}
		and
		\begin{align*}
			G_{zz} = & \left[\left(p'(V+r)-p'(V)-p''(V)r\right)V_z + \left(p'(V+r)-p'(V)\right)r_{z}\right]_z
			\\
			= & \left[p''(V+r)(V_z+r_{z}) - p''(V)V_z - p'''(V)V_z r - p''(V)r_{z}\right]V_z
			\\
			& + \left(p'(V+r)-p'(V)-p''(V)r\right)V_{zz} + \left[
			p''(V+r)(V_z+r_{z})
			\right.
			\\
			& \left. \quad -p''(V)V_z\right]r_{z} + \left(p'(V+r)-p'(V)\right)r_{zz}
			\\
			= & \left(p''(V+r)-p''(V)-p'''(V)r\right)V_z^2 + \left(p''(V+r)-p''(V)\right)r_{z}V_z
			\\
			& + \left(p'(V+r)-p'(V)-p''(V)r\right)V_{zz} + \left(p''(V+r)-p''(V)\right)V_z r_{z}
			\\
			& + p''(V+r)r_{z}^2 + \left(p'(V+r)-p'(V)\right)r_{zz}.
		\end{align*}
		Since $p \in C^{4}$, Taylor's expansion, together with \eqref{3-15}, gives the estimates of $|W_{zz}|$ and $|G_{zz}|$.
	\end{proof}
	
	\begin{lemma}\label{lemma3.3}
		Let the assumptions of Theorem \ref{theorem3.1} hold. If $N(T) \ll 1$, then there exits a constant $C>0$ independent of $T$ such that
		\begin{equation}
			\|\phi\|^2+\|\psi\|^2+\int_0^t \|r(\tau)\|_{w_1}^2 + \int_0^t \int d(V) \psi^2 \leq C (\|\phi_0\|^2 + \|\psi_0\|^2),
			\label{3-22}
		\end{equation}
		for any $t \in[0, T]$.
	\end{lemma}
	\begin{proof}
		Multiplying \eqref{3-6}$_{1}$ by $\phi$ and \eqref{3-6}$_{2}$ by $\frac{\psi}{-p'(V)}$, adding the resulting equations and integrating over $\mathbb{R} \times [0, t]$, we obtain
		\begin{equation}
			\begin{aligned}
				& \frac{1}{2} \int \left(\phi^2-\frac{\psi^2}{p'(V)}\right) + \int_0^t \int \frac{r^2}{V^{m+1}} + \frac{1}{2} \int_0^t \int d(V) \psi^2
				\\
				= & \frac{1}{2} \int \left(\phi_0^2-\frac{\psi_0^2}{p'(V)}\right) + \int_0^t \int \left(\frac{\phi r}{V^{m+1}} + \frac{s}{2} \phi^2 - \frac{s \psi^2}{2 p'(V)} + \phi \psi + W \phi\right)_z
				\\
				& - \int_0^t \int W r + \int_0^t \int G \frac{\psi}{p'(V)}.
			\end{aligned}
			\label{3-23}
		\end{equation}
		By \eqref{3-11}, the functions $\phi,\psi,r/V$ and $r/V^{m+1}$ belong to $H^1$, and hence tend to zero as $z \to \pm \infty$. Moreover, using $|W| \leq C r^{2}/V^{m+2} = C (r/V)(r/V^{m+1})$ and $1/V \leq C/V^{m+1}$, we deduce that $W \to 0$ as $z \to \pm \infty$. Consequently, 
		\begin{align}
			\int_0^t \int \left(\frac{\phi r}{V^{m+1}} + \frac{s}{2} \phi^2 - \frac{s \psi^2}{2 p'(V)} + \phi \psi + W \phi\right)_z = 0.
			\label{3-231}
		\end{align}
		Applying \eqref{3-15} and \eqref{3-18}, we obtain
		\begin{equation}
			\begin{aligned}
				\left| - \int_0^t \int W r + \int_0^t \int G \frac{\psi}{p'(V)} \right| \leq C \int_0^t \int \frac{|r|^3}{V^{m+2}} + C \int_0^t \int r^2 |\psi| \leq C N(t) \int_0^t \int \frac{r^2}{V}.
			\end{aligned}
			\label{3-24}
		\end{equation}
		Finally, substituting \eqref{3-231} and \eqref{3-24} into \eqref{3-23}, we arrive at
		\begin{align*}
			& \frac{1}{2} \int \left(\phi^2-\frac{\psi^2}{p'(V)}\right) + \int_0^t \int (1 - C N(t) V^{m}) \frac{r^2}{V^{m+1}} + \int_0^t \int \frac{s p''(V) V_z}{2 p'(V)^2} \psi^2
			\\
			\leq & \frac{1}{2} \int \left(\phi_0^2-\frac{\psi_0^2}{p'(V)}\right),
		\end{align*}
		which yields \eqref{3-22} for sufficiently small $N(t)$.
	\end{proof}
	
	\begin{lemma}\label{lemma3.4}
		Let the assumptions of Theorem \ref{theorem3.1} hold. If $N(T) \ll 1$, then there exits a constant $C>0$ independent of $T$ such that
		\begin{equation}
			\| r \|_{w_2}^2 + \int_0^t \|q\|_{w_1}^2 \leq C (\|\phi_0\|^2 + \|\psi_0\|^2 + \|r_{0}\|_{w_2}^2),
			\label{3-25}
		\end{equation}
		for any $t \in[0, T]$.
	\end{lemma}
	\begin{proof}
		To obtain the weighed space-time $L^{2}$ bound of $q$, we exploit the coupling between the two perturbation equations. Recalling \eqref{defs:r-q}, we have
		\begin{equation}
			r_{t} = r_{t} - s r_{z} + s r_{z} = q_{z} + s r_{z}.
			\label{3-251}
		\end{equation}
		Multiplying \eqref{3-251} by $\psi/V^{m+1}$ and \eqref{3-611}$_{2}$ by $r/V^{m+1}$, respectively, adding the resulting equations and integrating over $\mathbb{R}$, we obtain
		\begin{equation}
			\begin{aligned}
				\frac{\mathrm{d}}{\mathrm{d}t} \int \frac{r \psi}{V^{m+1}} = \int \left(\frac{r_{t} \psi}{V^{m+1}} + \frac{r \psi_{t}}{V^{m+1}}\right) = & \int \frac{q_{z} \psi}{V^{m+1}} + s \int \left(\frac{r}{V^{m+1}}\right)_{z} \psi + s \int \frac{(m+1)V_{z}}{V^{m+2}} r \psi
				\\
				& + s \int \frac{r \psi_{z}}{V^{m+1}} - \int \frac{r}{V^{m+1}} [p(V + r) - p(V)].
			\end{aligned}
			\label{3.3}
		\end{equation}
		By \eqref{3-11}, we know $\psi,r,r/V^{m+1},q/V^{m+1} \in H^{1}$, and hence
		\begin{equation}
			\psi,r,r/V^{m+1},q/V^{m+1} \to 0, \quad \hbox{as } z \to \pm \infty.
			\label{3.1211}
		\end{equation}
		Moreover, \eqref{3-151} and \eqref{3-6111} lead to $|\eta| \leq C|r|/V^{m+1}$, which, together with \eqref{3.1211}, implies
		\begin{equation}
			\eta \to 0, \quad \hbox{as } z \to \pm \infty.
			\label{3.121}
		\end{equation}
		Consequently, integration by parts gives
		\begin{equation}
			s \int \left(\frac{r}{V^{m+1}}\right)_{z} \psi + s \int \frac{r \psi_{z}}{V^{m+1}} = s \int \left(\frac{r \psi}{V^{m+1}}\right)_{z} = 0,
			\label{3.311}
		\end{equation}
		and
		\begin{equation}
			\begin{aligned}
				\int \frac{q_{z} \psi}{V^{m+1}} = & \int \left(\frac{q \psi}{V^{m+1}}\right)_{z} - \int \frac{q \psi_{z}}{V^{m+1}} + \int \frac{(m+1)V_{z}}{V^{m+2}} q \psi 
				\\
				= & - \int \frac{q^{2}}{V^{m+1}} + \int \frac{q \eta_{z}}{V^{m+1}} + \int \frac{(m+1)V_{z}}{V^{m+2}} q \psi,
			\end{aligned}
			\label{3.3111}
		\end{equation}
		where we have used $\psi_z = q - \eta_{z}$ in the last equality. Substituting \eqref{3.311} and \eqref{3.3111} ibto \eqref{3.3}, we arrive at
		\begin{equation}
			\begin{aligned}
				\frac{\mathrm{d}}{\mathrm{d}t} \int \frac{r \psi}{V^{m+1}} = & - \int \frac{q^{2}}{V^{m+1}} + \int \frac{q \eta_{z}}{V^{m+1}} + \int \frac{(m+1)V_{z}}{V^{m+2}} q \psi + s \int \frac{(m+1)V_{z}}{V^{m+2}} r \psi
				\\
				& - \int \frac{r}{V^{m+1}} [p(V + r) - p(V)].
			\end{aligned}
			\label{3.4}
		\end{equation}
		Although this identity provides the desired weighed space-time $L^{2}$ bound of $q$, the term $\int q \eta_{z}/V^{m+1}$ cannot be controlled by the estimates in Lemma \ref{lemma3.3}. Thus, we introduce an auxiliary functional. 
		
		To this end, we define
		\begin{equation*}
			\zeta := V^{m} \eta.
		\end{equation*}
		From the relationship between $\eta$ and $r$ in \eqref{3-6111}, we get 
		\begin{align}
			r = V [(-m V^{m} \eta + 1)^{-\frac{1}{m}} - 1],
			\label{3.711}
		\end{align}
		which implies
		\begin{equation}
			r_{t} = q_{z} + s r_{z} = q_{z} + sV_{z}[(-m V^{m} \eta + 1)^{-\frac{1}{m}} - 1] + sV(-m V^{m} \eta + 1)^{-\frac{1}{m}-1} (V^{m} \eta)_{z}.
			\label{3.9}
		\end{equation}
		Substituting \eqref{3.711} into \eqref{3.9}, we have
		\begin{equation}
			\begin{aligned}
				& V (-m \zeta + 1)^{-\frac{1}{m}-1} \zeta_{t} 
				\\
				= & q_{z} + sV_{z}[(-m \zeta + 1)^{-\frac{1}{m}} - 1] + sV(-m \zeta + 1)^{-\frac{1}{m}-1} \zeta_{z}.
			\end{aligned}
			\label{3.10}
		\end{equation}
		Multiplying \eqref{3.10} by $\zeta/V^{2m+1}$ and using \eqref{3.6}, we obtain
		\begin{equation}
			\begin{aligned}
				& \zeta (-m \zeta + 1)^{-\frac{1}{m}-1} \zeta_{t} /V^{2m}
				\\
				= & \zeta q_{z} /V^{2m+1} + s \zeta V_{z}[(-m \zeta + 1)^{-\frac{1}{m}} - 1] /V^{2m+1} + s \zeta (-m \zeta + 1)^{-\frac{1}{m}-1} \zeta_{z} /V^{2m}.
			\end{aligned}
			\label{3.101}
		\end{equation}
		Motivated by the structure of \eqref{3.101}, we introduce the auxiliary function
		\begin{align*}
			\Phi(\zeta) := \int_{0}^{\zeta} \theta (-m\theta+1)^{-\frac{1}{m}-1} \mathrm{d} \theta.
		\end{align*}
		It is obvious that
		\begin{equation}
			\Phi'(\zeta) = \zeta (-m \zeta + 1)^{-\frac{1}{m}-1},
			\label{3.6}
		\end{equation}
		and 
		\begin{equation*}
			\begin{aligned}
				\Phi''(\zeta) = (-m \zeta + 1)^{-\frac{1}{m}-1} + (m+1) (-m \zeta + 1)^{-\frac{1}{m}-2} \zeta = (\zeta+1) (-m \zeta + 1)^{-\frac{1}{m}-2}.
			\end{aligned}
		\end{equation*}
		Moreover, note that 
		\begin{equation*}
			\zeta =  - \frac{V^{m}}{m} \left[\frac{1}{(V+r)^{m}} - \frac{1}{V^{m}}\right] = - \frac{1}{m} \left[\frac{1}{(1+r/V)^{m}} - 1\right] = \int_{0}^{r/V} \frac{1}{(1+\theta)^{m+1}} \mathrm{d}\theta.
		\end{equation*}
		Since $|r/V| \leq C N(t) \ll 1$, we know that there exist positive constants $c$ and $C$, independent of $z$ and $t$, such that
		\begin{equation}
			c \left|\frac{r}{V}\right| \leq |\zeta| \leq C \left|\frac{r}{V}\right|.
			\label{3-55}
		\end{equation}
		Since $\Phi(0) = \Phi'(0) = 0$ and $\Phi''(0) \neq 0$, Taylor's theorem yields
		\begin{equation}
			c \left(\frac{r}{V}\right)^{2} \leq \Phi(\zeta) \leq C \left(\frac{r}{V}\right)^{2}.
			\label{3.33}
		\end{equation}
		Therefore, 
		\begin{equation}
			c \frac{r^{2}}{V^{2(m+1)}} \leq \frac{\Phi(\zeta)}{V^{2m}} \leq C \frac{r^{2}}{V^{2(m+1)}}.
			\label{3.331}
		\end{equation}
		Then, we derive an estimate for $\Phi(\zeta)/V^{2m}$. Using \eqref{3.101} and $\zeta = V^{m} \eta$, we obtain 
		\begin{equation}
			\begin{aligned}
				\frac{\mathrm{d}}{\mathrm{d}t} \int \frac{\Phi(\zeta)}{V^{2m}} = & \int \frac{\zeta}{V^{2m}} (-m \zeta + 1)^{-\frac{1}{m}-1} \zeta_{t} 
				\\
				= & \int \frac{\zeta  q_{z}}{V^{2m+1}} + \int \frac{sV_{z} \zeta}{V^{2m+1}} [(-m \zeta + 1)^{-\frac{1}{m}} - 1] + \int \frac{s \zeta}{V^{2m}}(-m \zeta + 1)^{-\frac{1}{m}-1} \zeta_{z}
				\\
				= & \int \frac{\eta  q_{z}}{V^{m+1}} + \int \frac{sV_{z} \eta r}{V^{m+2}} + \int \frac{s \zeta}{V^{2m}}(-m \zeta + 1)^{-\frac{1}{m}-1} \zeta_{z}.
			\end{aligned}
			\label{3.12}
		\end{equation}
		Integration by parts, together with \eqref{3.1211} and \eqref{3.121}, gives
		\begin{equation*}
			\int \frac{\eta  q_{z}}{V^{m+1}} = \int \left(\frac{\eta q}{V^{m+1}}\right)_{z} - \int \frac{\eta_{z} q}{V^{m+1}} + \int \frac{(m+1)V_{z} \eta q}{V^{m+2}} = - \int \frac{\eta_{z} q}{V^{m+1}} + \int \frac{(m+1)V_{z} \eta q}{V^{m+2}}. 
		\end{equation*}
		Furthermore, \eqref{3.1211} and \eqref{3.331} imply $\Phi(\zeta)/V^{2m} \to 0$ as $z \to \pm \infty$. This along with $\Phi_z(\zeta) = \Phi'(\zeta) \zeta_{z} = \zeta (-m \zeta + 1)^{-\frac{1}{m}-1} \zeta_{z}$ gives
		\begin{align*}
			\int \frac{s \zeta}{V^{2m}}(-m \zeta + 1)^{-\frac{1}{m}-1} \zeta_{z} = s \int \frac{\Phi_{z}(\zeta)}{V^{2m}} = s \int \left(\frac{\Phi(\zeta)}{V^{2m}}\right)_{z} + 2 s m \int \frac{V_{z} \Phi(\zeta)}{V^{2m+1}} = 2 s m \int \frac{V_{z} \Phi(\zeta)}{V^{2m+1}}.
		\end{align*}
		Substituting the above identities into \eqref{3.12}, we have
		\begin{equation}
			\begin{aligned}
				\frac{\mathrm{d}}{\mathrm{d}t} \int \frac{\Phi(\zeta)}{V^{2m}} = & - \int \frac{\eta_{z} q}{V^{m+1}} + \int \frac{(m+1)V_{z} \eta q}{V^{m+2}} + \int \frac{sV_{z} \eta r}{V^{m+2}} + 2 s m \int \frac{V_{z} \Phi(\zeta)}{V^{2m+1}}.
			\end{aligned}
			\label{3.15}
		\end{equation}
		Adding \eqref{3.4} and \eqref{3.15}, we deduce that
		\begin{align*}
			\frac{\mathrm{d}}{\mathrm{d}t} \int \frac{\Phi(\zeta)}{V^{2m}} + \frac{\mathrm{d}}{\mathrm{d}t} \int \frac{r \psi}{V^{m+1}}  + \int \frac{q^{2}}{V^{m+1}} = & \int \frac{(m+1)V_{z}}{V^{m+2}} q (\psi + \eta) - \int \frac{r}{V^{m+1}} [p(V + r) - p(V)]  
			\\
			& + s \int \frac{(m+1)V_{z}}{V^{m+2}} r \psi + \int \frac{sV_{z} \eta r}{V^{m+2}} + 2 s m \int \frac{V_{z} \Phi(\zeta)}{V^{2m+1}}.
		\end{align*}
		Since $q = (\psi + \eta)_{z}$, integration by parts yields
		\begin{align*}
			\int \frac{(m+1)V_{z}}{V^{m+2}} q (\psi + \eta) = & \frac{1}{2} \int \left[\frac{(m+1)V_{z}}{V^{m+2}} (\psi + \eta)^{2}\right]_{z} - \frac{1}{2} \int \left(\frac{(m+1)V_{z}}{V^{m+2}}\right)_{z} (\psi + \eta)^{2} 
			\\
			= & - \frac{1}{2} \int \left(\frac{(m+1)V_{z}}{V^{m+2}}\right)_{z} (\psi + \eta)^{2}.
		\end{align*}
		Combining the above identity, we conclude that
		\begin{equation}
			\begin{aligned}
				& \frac{\mathrm{d}}{\mathrm{d}t} \int \frac{\Phi(\zeta)}{V^{2m}} + \frac{\mathrm{d}}{\mathrm{d}t} \int \frac{r \psi}{V^{m+1}} + \int \frac{q^{2}}{V^{m+1}} 
				\\
				= & - \frac{1}{2} \int \left(\frac{(m+1)V_{z}}{V^{m+2}}\right)_{z} (\psi + \eta)^{2} - \int \frac{r}{V^{m+1}} [p(V + r) - p(V)] + s \int \frac{(m+1)V_{z}}{V^{m+2}} r \psi
				\\
				& + \int \frac{sV_{z} \eta r}{V^{m+2}} + 2 s m \int \frac{V_{z} \Phi(\zeta)}{V^{2m+1}}.
			\end{aligned}
			\label{3.18}
		\end{equation}
		Next, we need to eliminate the unfavorable term $s \int \frac{(m+1)V_{z}}{V^{m+2}} r \psi$. Multiplying \eqref{3-6}$_{2}$ by $\frac{s(m+1)V_{z}}{p'(V)V^{m+2}} \psi$ gives
		\begin{equation}
			\begin{aligned}
				\frac{s}{2} \frac{\mathrm{d}}{\mathrm{d}t} \int \frac{(m+1)V_{z}}{p'(V) V^{m+2}} \psi^{2} = \int\frac{s^{2} (m+1)V_{z}}{p'(V) V^{m+2}} \psi \psi_{z} - s \int \frac{(m+1)V_{z}}{V^{m+2}} \psi r - s \int \frac{(m+1)V_{z}}{p'(V) V^{m+2}} \psi G.
			\end{aligned}
			\label{3.20}
		\end{equation}
		Noting \eqref{3.1211} and using integration by parts, we get
		\begin{equation*}
			\int\frac{s^{2} (m+1)V_{z}}{p'(V) V^{m+2}} \psi \psi_{z} = \int \left(\frac{s^{2}(m+1)V_{z} \psi^{2}}{2 p'(V) V^{m+2}}\right)_{z} - \int \left(\frac{s^{2}(m+1)V_{z}}{2 p'(V) V^{m+2}}\right)_{z} \psi^{2} = - \int \left(\frac{s^{2}(m+1)V_{z}}{2 p'(V) V^{m+2}}\right)_{z} \psi^{2}.
		\end{equation*}
		Combining the above equation, \eqref{3.18} and \eqref{3.20}, we deduce that
		\begin{equation}
			\begin{aligned}
				& \frac{\mathrm{d}}{\mathrm{d}t} \int \frac{\Phi(\zeta)}{V^{2m}} + \frac{\mathrm{d}}{\mathrm{d}t} \int \frac{r \psi}{V^{m+1}} + \frac{s}{2} \frac{\mathrm{d}}{\mathrm{d}t} \int \frac{(m+1)V_{z}}{p'(V) V^{m+2}} \psi^{2} + \int \frac{q^{2}}{V^{m+1}} 
				\\
				= & - \frac{1}{2} \int \left(\frac{(m+1)V_{z}}{V^{m+2}}\right)_{z} (\psi + \eta)^{2} - \int \frac{r}{V^{m+1}} [p(V + r) - p(V)] + \int \frac{sV_{z} \eta r}{V^{m+2}} 
				\\
				& + 2 s m \int \frac{V_{z} \Phi(\zeta)}{V^{2m+1}} - \int \left(\frac{s^{2}(m+1)V_{z}}{2 p'(V) V^{m+2}}\right)_{z} \psi^{2} - s \int \frac{(m+1)V_{z}}{p'(V) V^{m+2}} \psi G
				\\
				=: & \sum_{i=1}^{6} I_{i}.
			\end{aligned}
			\label{3.21}
		\end{equation}
		It remains to estimate the six terms $\{I_{i}\}_{i=1}^{6}$ on the right-hand side of \eqref{3.21}. By Taylor's expansion, it holds that
		\begin{equation*}
			\left| \eta \right| = \left| -\frac{1}{m}\left[\frac{1}{(V+r)^{m}} - \frac{1}{V^{m}}\right] \right| \leq C \left|\frac{r}{V^{m+1}}\right|.
		\end{equation*}
		This along with \eqref{3-16} and \eqref{profile:coeff} leads to
		\begin{equation*}
			\begin{aligned}
				|I_{1} + I_{5}| \leq C \int |V_{z}| [(\psi + \eta)^{2} + \psi^{2}] \leq C \int \frac{r^{2}}{V^{m}} + C \int |V_{z}| \psi^{2},
			\end{aligned}
		\end{equation*}
		and 
		\begin{equation}
			|I_{3}| \leq C \int \frac{r^{2}}{V^{m+1}}.
			\label{3.25}
		\end{equation}
		For $I_{2} + I_{4} + I_{6}$, using \eqref{3-16}, \eqref{3-18}, Taylor's expansion and the fact $\| \psi \|_{L^{\infty}} \leq N(t)$, we obtain
		\begin{equation}
			|I_{2} + I_{4} + I_{6}| \leq C \int \frac{r^{2}}{V^{m+1}} + C N(t) \int r^{2}.
			\label{3.27}
		\end{equation}
		Substituting \eqref{3.25}--\eqref{3.27} into \eqref{3.21}, from $1/V^{m} \leq C/V^{m+1}$, we infer that
		\begin{equation}
			\begin{aligned}
				& \frac{\mathrm{d}}{\mathrm{d}t} \int \frac{\Phi(\zeta)}{V^{2m}} + \frac{\mathrm{d}}{\mathrm{d}t} \int \frac{r \psi}{V^{m+1}} - \frac{s}{2} \frac{\mathrm{d}}{\mathrm{d}t} \int \frac{(m+1)V_{z}}{-p'(V) V^{m+2}} \psi^{2} + \int \frac{q^{2}}{V^{m+1}} 
				\\
				\leq & C \int \frac{r^{2}}{V^{m+1}} + C \int r^{2} + C \int |V_{z}| \psi^{2}.
			\end{aligned}
			\label{3.28}
		\end{equation}
		Integrating \eqref{3.28} over $[0,t]$ and using \eqref{3.331}, we get
		\begin{equation}
			\begin{aligned}
				& c \int \frac{r^{2}}{V^{2(m+1)}} + \int \frac{r \psi}{V^{m+1}} + \frac{s}{2} \int \frac{(m+1)V_{z}}{p'(V) V^{m+2}} \psi^{2} + \int_{0}^{t} \int \frac{q^{2}}{V^{m+1}}
				\\
				\leq & C \int \left(\frac{r_{0}^{2}}{V^{2(m+1)}} + \frac{r_{0} \psi_{0}}{V^{m+1}}\right) + \frac{s}{2} \int \frac{(m+1)V_{z}}{p'(V) V^{m+2}} \psi_{0}^{2} + C \int_{0}^{t} \int \frac{r^{2}}{V^{m+1}} + C \int_{0}^{t} \int |V_{z}| \psi^{2}.
			\end{aligned}
			\label{3.29}
		\end{equation}
		Since $s < 0$, $V_{z} < 0$ and $-p'(V) > 0$, it follows that $\int \frac{s(m+1)V_{z}}{2 p'(V) V^{m+2}} \psi^{2} < 0$. By $p \in C^{4}([0,\infty)$, we have
		\begin{equation}
			\left|\int \frac{s(m+1)V_{z}}{2 p'(V) V^{m+2}} \psi^{2}\right| \leq C \int \psi^{2}.
			\label{3.291}
		\end{equation}
		In addition, Young's inequality gives
		\begin{equation}
			\left| \int \frac{r \psi}{V^{m+1}} \right| \leq \frac{c}{2} \int \left(\frac{r}{V^{m+1}}\right)^{2} + C \int \psi^{2}.
			\label{3.35}
		\end{equation}
		Applying \eqref{3.291}, \eqref{3.35} and Lemma \ref{lemma3.3}, \eqref{3.29} becomes
		\begin{align*}
			\int \left(\frac{r}{V^{m+1}}\right)^{2} + \int_{0}^{t} \int \frac{q^{2}}{V^{m+1}} \leq & C (\|\phi_0\|^2 + \|\psi_0\|^2 + \|r_{0}\|_{w_2}^2) + C \int \psi^{2} + C \int_{0}^{t} \int \left(\frac{r^{2}}{V^{m+1}} + |V_{z}| \psi^{2}\right)
			\nonumber
			\\
			\leq & C (\|\phi_0\|^2 + \|\psi_0\|^2 + \|r_{0}\|_{w_2}^2),
		\end{align*}
		which completes the proof. 
	\end{proof}

	\subsection{Singular weighted energy estimates}
	Differentiating \eqref{3-61} with respect to $t$ and $z$, respectively, we obtain
	\begin{equation*}
		q_{t}=\psi_{zt}+\eta_{zt}, \qquad sq_{z} = s \psi_{zz} + s \eta_{zz}.
	\end{equation*}
	Consequently,
	\begin{equation}
		q_t = sq_z + (\psi_{t} - s\psi_{z})_{z} + (\eta_{t} - s\eta_{z})_{z}.
		\label{4.71}
	\end{equation}
	By \eqref{3-611}$_{2}$, we have $(\psi_{t} - s\psi_{z})_{z} = - [p(V+r)-p(V)]_{z}$. From \eqref{3-6111} and $r_t=s r_z+q_z$, we find
	\begin{equation*}
		\begin{aligned}
			\eta_t-s\eta_z=\frac{r_{t}}{(V+r)^{m+1}}-s\left(\frac{V_z+r_{z}}{(V+r)^{m+1}}-\frac{V_z}{V^{m+1}}\right)=\frac{q_{z}}{(V+r)^{m+1}} -s\left(\frac{V_z}{(V+r)^{m+1}}-\frac{V_z}{V^{m+1}}\right).
		\end{aligned}
	\end{equation*}
	We then rewrite \eqref{4.71} as
	\begin{equation}
		q_t = sq_z + \left[\frac{q_{z}}{(V+r)^{m+1}} -s\left(\frac{V_z}{(V+r)^{m+1}}-\frac{V_z}{V^{m+1}}\right) - (p(V+r)-p(V))\right]_z.
		\label{4.7}
	\end{equation}
	Moreover, since $(V+r)^{-(m+1)} = V^{-(m+1)} V^{m+1}/(V+r)^{m+1} = V^{-(m+1)} (1+r/V)^{-(m+1)}$, the smallness condition \eqref{3-151} ensures that there exist positive constant $c_{1}$ and $C$, independent of $z$ and $t$, such that
	\begin{equation}
		\frac{c_{1}}{V^{m+1}} \leq \frac{1}{(V+r)^{m+1}} \leq \frac{C}{V^{m+1}}.
		\label{4.711}
	\end{equation}
	To handle the singular coefficients at $z = +\infty$, we introduce a smooth cut-off function $\kappa \in C_c^\infty(\R)$ satisfying
	\begin{equation}
		\begin{cases}
			\kappa(z) = 1, \, \;\hbox{for } |z| < L; \quad \kappa(z) = 0, \, \;\hbox{for } |z| > 2L;
			\\
			0 \leq \kappa(z) \leq 1, \quad |\kappa_{z}(z)| \leq \frac{C}{L}, \,\; \hbox{for } z \in \mathbb{R},
		\end{cases}
		\label{4.10}
	\end{equation}
	where $L$ is sufficiently large. Before proceeding with the weighted energy estimates, we first establish the time regularity required for the subsequent estimates. For any fixed $L>0$, let $I_L:=(-3L,3L)$. In view of \eqref{3-16}, \eqref{3-151} and \eqref{3-11}, the coefficients in \eqref{4.7} are bounded on $I_L\times(0,T)$. Furthermore, we get from the mean value theorem that
	\begin{equation*}
		\left|s\left(\frac{V_z}{(V+r)^{m+1}}-\frac{V_z}{V^{m+1}}\right) + (p(V+r)-p(V))\right| \leq C|r|.
	\end{equation*}
	Since $r,q\in C([0,T];H^1_{w_2}(\mathbb R))$, we deduce from \eqref{4.7} that
	\begin{equation*}
		q_t\in L^2(0,T;H^{-1}(I_L)).
	\end{equation*}
	Consequently, for each fixed \(L>0\), the Lions--Magenes chain rule justifies the time differentiation of the weighted $L^2$-norm of $q$ with weight $\kappa^2/V^{2(m+1)}$. Thus, all subsequent time-dependent energy calculations are justified in the corresponding duality sense.
	
	\begin{lemma}\label{lemma3.6}
		Let the assumptions of Theorem \ref{theorem3.1} hold. If $N(T) \ll 1$, then there exists a constant $C>0$ independent of $T$ such that
		\begin{equation}
			\|q(t)\|_{w_2}^2+\int_0^t\|q_z\|_{w_3}^2\,\mathrm{d}\tau
			\leq C\bigl(\|\phi_0\|^2+\|\psi_0\|^2
			+\|r_0\|_{w_2}^2+\|q_0\|_{w_2}^2\bigr),
			\label{3-46}
		\end{equation}
		for any $t \in[0, T]$.
	\end{lemma}
	\begin{proof}
		Multiplying both sides of \eqref{4.7} by $\frac{q}{V^{2(m+1)}}$ gives
		\begin{equation}
			\begin{aligned}
				\frac{1}{2} \frac{\mathrm{d}}{\mathrm{d}t} \frac{q^{2}}{V^{2(m+1)}} = & \left[\frac{q_{z} q}{(V+r)^{m+1}V^{2(m+1)}} -s\left(\frac{V_z}{(V+r)^{m+1}}-\frac{V_z}{V^{m+1}}\right)   \frac{q}{V^{2(m+1)}} \right. 
				\\
				& \left. - (p(V+r)-p(V)) \frac{q}{V^{2(m+1)}} + \frac{s}{2} \frac{q^{2}}{V^{2(m+1)}}\right]_{z} + s \frac{(m+1)V_{z} q^{2}}{V^{2m+3}}
				\\
				& - \frac{q_{z}}{(V+r)^{m+1}} \left(\frac{q_{z}}{V^{2(m+1)}} - \frac{2 (m+1) V_{z} q}{V^{2m+3}}\right) + \bigg[s\left(\frac{V_z}{(V+r)^{m+1}}-\frac{V_z}{V^{m+1}}\right) 
				\\
				& + (p(V+r)-p(V))\bigg] \left(\frac{q_{z}}{V^{2(m+1)}} - \frac{2 (m+1) V_{z} q}{V^{2m+3}}\right).
			\end{aligned}
			\label{3-4611}
		\end{equation}
		Next, multiplying the above identity by $\kappa^{2}(z)$ defined in \eqref{4.10} and integrating the result over $\mathbb{R} \times [0,t]$, we obtain
		\begin{align*}
			& \frac{1}{2} \int_{-2L}^{2L} \frac{\kappa^{2} q^{2}}{V^{2(m+1)}} + \int_{0}^{t} \int_{-2L}^{2L} \frac{\kappa^{2} q_{z}^{2}}{V^{2(m+1)}(V+r)^{m+1}}
			\\
			= & \frac{1}{2} \int_{-2L}^{2L} \frac{\kappa^{2} q_{0}^{2}}{V^{2(m+1)}} + s \int_{0}^{t} \int_{-2L}^{2L} \frac{\kappa^{2} (m+1)V_{z} q^{2}}{V^{2m+3}} + 2 \int_{0}^{t} \int_{-2L}^{2L} \kappa^{2} \frac{(m+1) V_{z} q q_{z}}{V^{2m+3}(V+r)^{m+1}} 
			\\
			& + \int_{0}^{t} \int_{-2L}^{2L} \kappa^{2} \bigg[s\left(\frac{V_z}{(V+r)^{m+1}}-\frac{V_z}{V^{m+1}}\right) + (p(V+r)-p(V))\bigg] \left(\frac{q_{z}}{V^{2(m+1)}} - \frac{2 (m+1) V_{z} q}{V^{2m+3}}\right)  
			\\
			& - \int_{0}^{t} \int_{-2L}^{2L} 2 \kappa \kappa_{z} \left[\frac{q_{z} q}{(V+r)^{m+1}V^{2(m+1)}} -\left(\frac{s V_z}{(V+r)^{m+1}}-\frac{s V_z}{V^{m+1}} + (p(V+r)-p(V))\right) \right.
			\\
			& \left. \times \frac{q}{V^{2(m+1)}} + \frac{s}{2} \frac{q^{2}}{V^{2(m+1)}}\right].
		\end{align*}
		We next estimate the last term of the above equation. Using \eqref{3-16}, \eqref{4.711}, Taylor's expansion, Young's inequality and the bound $|\kappa_{z}| \leq C/L$, we deduce that
		\begin{align*}
			& \bigg| \int_{0}^{t} \int_{-2L}^{2L} 2 \kappa \kappa_{z} \bigg[ \frac{q_{z} q}{(V+r)^{m+1}V^{2(m+1)}} -\bigg(\frac{s V_z}{(V+r)^{m+1}}-\frac{s V_z}{V^{m+1}} 
			\\
			& + (p(V+r)-p(V))\bigg) \frac{q}{V^{2(m+1)}} + \frac{s}{2} \frac{q^{2}}{V^{2(m+1)}} \bigg] \bigg|
			\\
			\leq & \frac{C}{L} \int_{0}^{t} \int_{-2L}^{2L} \kappa \left( \frac{|q| |q_{z}|}{V^{3(m+1)}} +  \frac{|V_{z}||r| |q|}{V^{3m+4}} + \frac{|r| |q|}{V^{2(m+1)}} + \frac{q^{2}}{V^{2(m+1)}} \right)
			\\
			\leq & \frac{C}{\sqrt{L}} \int_0^t \int_{-2L}^{2L} \kappa^{2} \frac{q_z^2}{V^{3(m+1)}}+\frac{C}{L^{\frac{3}{2}}} \int_{0}^{t} \int_{-2L}^{2L} \frac{q^{2}}{V^{3(m+1)}}+\frac{C}{L}\int_{0}^{t} \int_{-2L}^{2L} \frac{q^{2}+r^{2}}{V^{2(m+1)}}.
		\end{align*}
		By Theorem \ref{theorem1.1}, it is clear that $V(z) \sim (a_{+}z)^{-\frac{1}{m+1}}$ as $z \to +\infty$, which implies $1/V^{m+1} \leq CL$ for $z \in [-2L,2L]$. This, together with $\|r\|_{w_{2}}^2,\|q\|_{w_{2}}^2 \leq N^2(t)$ and Lemma \ref{lemma3.3}, allows the inequality above to be estimated as
		\begin{equation}
			\begin{aligned}
				& \bigg| \int_{0}^{t} \int_{-2L}^{2L} 2 \kappa \kappa_{z} \bigg[ \frac{q_{z} q}{(V+r)^{m+1}V^{2(m+1)}} -\bigg(\frac{s V_z}{(V+r)^{m+1}}-\frac{s V_z}{V^{m+1}} 
				\\
				& + (p(V+r)-p(V))\bigg) \frac{q}{V^{2(m+1)}} + \frac{s}{2} \frac{q^{2}}{V^{2(m+1)}} \bigg] \bigg|
				\\
				\leq & \frac{C}{\sqrt{L}} \int_0^t \int_{-2L}^{2L} \kappa^{2} \frac{q_z^2}{V^{3(m+1)}} + \frac{C T N^{2}(T)}{\sqrt{L}}.
			\end{aligned}
			\label{3-371}
		\end{equation}
		It follows from \eqref{4.711} and \eqref{3-371} that
		\begin{align*}
			& \frac{1}{2} \int_{-2L}^{2L} \kappa^{2} \frac{q^{2}}{V^{2(m+1)}} + \left(c_{1} - \frac{C}{\sqrt{L}}\right) \int_{0}^{t} \int_{-2L}^{2L} \kappa^{2} \frac{q_{z}^{2}}{V^{3(m+1)}}
			\\
			\leq & \frac{1}{2} \int_{-2L}^{2L} \kappa^{2} \frac{q_{0}^{2}}{V^{2(m+1)}} + \frac{C T N^{2}(T)}{\sqrt{L}} + s \int_{0}^{t} \int_{-2L}^{2L} \kappa^{2} \frac{(m+1)V_{z} q^{2}}{V^{2m+3}} + 2 \int_{0}^{t} \int_{-2L}^{2L} \kappa^{2} \frac{(m+1) V_{z} q q_{z}}{V^{2m+3}(V+r)^{m+1}} 
			\\
			& + \int_{0}^{t} \int_{-2L}^{2L} \kappa^{2} \bigg[s\left(\frac{V_z}{(V+r)^{m+1}}-\frac{V_z}{V^{m+1}}\right) + (p(V+r)-p(V))\bigg] \left(\frac{q_{z}}{V^{2(m+1)}} - \frac{2 (m+1) V_{z} q}{V^{2m+3}}\right).
		\end{align*}
		Passing the limit $L \to +\infty$ gives
		\begin{equation*}
			\begin{aligned}
				& \frac{1}{2} \int \frac{q^{2}}{V^{2(m+1)}} + \frac{3 c_{1}}{4} \int_{0}^{t} \int \frac{q_{z}^{2}}{V^{3(m+1)}}
				\\
				\leq & \frac{1}{2} \int \frac{q_{0}^{2}}{V^{2(m+1)}} + \int_{0}^{t} \int \left| \frac{s (m+1)V_{z} q^{2}}{V^{2m+3}} \right| + 2 \int_{0}^{t} \int \left| \frac{(m+1) V_{z} q q_{z}}{V^{2m+3}(V+r)^{m+1}} \right|
				\\
				& + \int_{0}^{t} \int \left|\bigg[s\left(\frac{V_z}{(V+r)^{m+1}}-\frac{V_z}{V^{m+1}}\right) + (p(V+r)-p(V))\bigg] \left(\frac{q_{z}}{V^{2(m+1)}} - \frac{2 (m+1) V_{z} q}{V^{2m+3}}\right)\right|.
			\end{aligned}
		\end{equation*}
		Thanks to $|V_{z}| \leq C V^{m+2}$ and $1/(V+r)^{m+1} \leq C/V^{m+1}$ in \eqref{4.711}, using Taylor's expansion and Young's inequality, one can see that
		\begin{equation*}
			\begin{aligned}
				& \int_{0}^{t} \int \left|\frac{s (m+1)V_{z} q^{2}}{V^{2m+3}}\right| + 2 \int_{0}^{t} \int \left|\frac{(m+1) V_{z} q q_{z}}{V^{2m+3}(V+r)^{m+1}}\right| 
				\\
				& + \int_{0}^{t} \int \left|\bigg[s\left(\frac{V_z}{(V+r)^{m+1}}-\frac{V_z}{V^{m+1}}\right) + (p(V+r)-p(V))\bigg] \left(\frac{q_{z}}{V^{2(m+1)}} - \frac{2 (m+1) V_{z} q}{V^{2m+3}}\right)\right|
				\\
				\leq & C \int_{0}^{t} \int \frac{q^{2}}{V^{m+1}} + C \int_{0}^{t} \int \frac{|q| |q_z|}{V^{2(m+1)}} + C \int_{0}^{t} \int \frac{|r| |q_z|}{V^{2(m+1)}} + C \int_{0}^{t} \int \frac{|r| |q|}{V^{m+1}}
				\\
				\leq & \frac{c_{1}}{4} \int_{0}^{t} \int \frac{q_z^2}{V^{3(m+1)}} + C \int_{0}^{t} \int \frac{q^2+r^2}{V^{m+1}}.
			\end{aligned}
		\end{equation*}
		In view of the above estimates, Lemmas \ref{lemma3.3} and \ref{lemma3.4}, it follows that
		\begin{equation}
			\begin{aligned}
				\frac{1}{2} \int \frac{q^{2}}{V^{2(m+1)}} + \frac{c_{1}}{2} \int_{0}^{t} \int \frac{q_{z}^{2}}{V^{3(m+1)}} \leq & C \int \frac{q_{0}^{2}}{V^{2(m+1)}} + C \int_{0}^{t} \int \frac{q^{2}}{V^{m+1}} + C \int_{0}^{t} \int \frac{r^{2}}{V^{m+1}}
				\\
				\leq & C (\|\phi_0\|^2 + \|\psi_0\|^2 + \|r_{0}\|_{w_2}^2 +  \|q_{0}\|_{w_{2}}^2).
			\end{aligned}
			\label{3-47}
		\end{equation}
		The proof of this lemma is complete.
	\end{proof}
	
	\begin{lemma}\label{lemma3.5}
		Let the assumptions of Theorem \ref{theorem3.1} hold. If $N(T) \ll 1$, then there exists a constant $C>0$ independent of $T$ such that
		\begin{equation}
			\|r(t)\|^2+\|\psi_z(t)\|^2
			+\int_0^t \|r_z\|_{w_1}^2 \mathrm{d}\tau + \int_0^t \int d(V)\psi_z^2 \mathrm{d}\tau
			\leq C(\|\phi_0\|_1^2+\|\psi_0\|_1^2),
			\label{3-41}
		\end{equation}
		for any $t \in[0, T]$.
	\end{lemma}
	\begin{proof}
		Differentiating \eqref{3-6} with respect to $z$, we get
		\begin{equation}
			\begin{cases}
				r_t = sr_{z} + \psi_{zz} + (r/V^{m+1})_{zz} + W_{zz},
				\\
				\psi_{zt} = s\psi_{zz} - p''(V) V_{z}r - p'(V)r_{z} - G_{z}.
			\end{cases}
			\label{3-42}
		\end{equation}
		Multiplying the first equation of \eqref{3-42} by $r$ and the second equation by $(-\frac{\psi_z}{p'(V)})$, and then adding them, we derive
		\begin{equation}
			\begin{aligned}
				& \frac{1}{2} \frac{\mathrm{d}}{\mathrm{d}t} \left(r^2 - \frac{\psi_z^2}{p'(V)}\right) 
				\\
				= & \left[ \left(\frac{r_{z}}{V^{m+1}} - \frac{(m+1)V_{z} r}{V^{m+2}}\right) r + \frac{s}{2} r^{2} + \psi_{z} r + W_{z} r + \frac{s\psi_z^2}{-2p'(V)} \right]_{z} - \frac{r_{z}^2}{V^{m+1}}
				\\
				& + \frac{(m+1) V_z}{V^{m+2}} r r_{z} - \frac{s p''(V) V_z}{2 p'(V)^2} \psi_z^2 - W_{z} r_{z} + \frac{p''(V) V_{z}}{p'(V)} r \psi_z +  G_z \frac{\psi_z}{p'(V)}.
			\end{aligned}
			\label{3-4211}
		\end{equation}
		We multiply \eqref{3-4211} by the cut-off function $\kappa^{2}(z)$ and integrate over $\mathbb{R} \times [0,t]$. Then, we derive
		\begin{equation}
			\begin{aligned}
				& \frac{1}{2} \int_{-2L}^{2L} \kappa^{2} \left(r^2 - \frac{\psi_z^2}{p'(V)}\right) + \int_{0}^{t} \int_{-2L}^{2L} \kappa^{2} \frac{r_{z}^2}{V^{m+1}} + \int_{0}^{t} \int_{-2L}^{2L} \kappa^{2} \frac{s p''(V) V_z}{2 p'(V)^2} \psi_z^2
				\\
				= & \frac{1}{2} \int_{-2L}^{2L} \kappa^{2} \left(r_{0}^2 - \frac{\psi_{0z}^2}{p'(V)}\right) + \int_{0}^{t} \int_{-2L}^{2L} \kappa^{2} \left(\frac{(m+1) V_z}{V^{m+2}} r r_{z} - W_{z} r_{z} + \frac{p''(V) V_{z}}{p'(V)} r \psi_z \right.
				\\
				& \left. +  G_z \frac{\psi_z}{p'(V)}\right) - \int_{0}^{t} \int_{-2L}^{2L} 2 \kappa \kappa_{z} \bigg[\left(\frac{r_{z}}{V^{m+1}} - \frac{(m+1)V_{z} r}{V^{m+2}}\right) r + \frac{s}{2} r^{2} + \psi_{z} r + W_{z} r - \frac{s\psi_z^2}{2p'(V)}\bigg].
			\end{aligned}
			\label{3-4311}
		\end{equation}
		We next estimate the last term of \eqref{3-4311}. Combining \eqref{3-18} and Lemma \ref{lemma3.3} and utilizing the estimates $|\kappa_{z}| \leq C/L$, $\|\psi_{z}\|, \|r/V^{m+1}\|_{L^{\infty}}, \|r_{z}/V^{m+1}\|_{L^{\infty}} \leq C N(t)$,  $|V_{z}| \leq C V^{m+2}$, one gets
		\begin{equation*}
			\begin{aligned}
				& \left|\int_{0}^{t} \int_{-2L}^{2L} 2 \kappa \kappa_{z} \bigg[\left(\frac{r_{z}}{V^{m+1}} - \frac{(m+1)V_{z} r}{V^{m+2}}\right) r + \frac{s}{2} r^{2} + \psi_{z} r + W_{z} r + \frac{s\psi_z^2}{-2p'(V)}\bigg]\right|
				\\
				\leq & \frac{C}{L} \int_{0}^{t} \int_{-2L}^{2L} \kappa \bigg(\frac{|r_{z}| |r|}{V^{m+1}} + r^{2} + |\psi_{z}| |r| + \frac{|r|^{3}}{V} + \frac{|r_{z}| |r|^{2}}{V^{m+2}} + \psi_{z}^{2} \bigg)
				\\
				\leq & \frac{C(1+N(T))}{L} \int_{0}^{t} \int_{-2L}^{2L} \kappa^{2} \frac{r_{z}^{2}}{V^{m+1}} + \frac{C(1+N(T))}{L} \int_{0}^{t} \int_{-2L}^{2L} \frac{r^{2}}{V^{m+1}} + \frac{C}{L} \int_{0}^{t} \int_{-2L}^{2L} \psi_{z}^{2}
				\\
				\leq & \frac{C(1+N(T))}{L} \int_{0}^{t} \int_{-2L}^{2L} \kappa^{2} \frac{r_{z}^{2}}{V^{m+1}} + \frac{C(1+N(T))}{L} (\|\phi_0\|^2 + \|\psi_0\|^2) + \frac{C T N^{2}(T)}{L}.
			\end{aligned}
		\end{equation*}
		Therefore, \eqref{3-4311} is estimated as
		\begin{equation*}
			\begin{aligned}
				& \frac{1}{2} \int_{-2L}^{2L} \kappa^{2} \left(r^2 - \frac{\psi_z^2}{p'(V)}\right) + \left(1 - \frac{C(1+N(T))}{L}\right) \int_{0}^{t} \int_{-2L}^{2L} \kappa^{2} \frac{r_{z}^2}{V^{m+1}} + \int_{0}^{t} \int_{-2L}^{2L} \kappa^{2} \frac{s p''(V) V_z}{2 p'(V)^2} \psi_z^2
				\\
				\leq & \frac{1}{2} \int_{-2L}^{2L} \kappa^{2} \left(r_{0}^2 - \frac{\psi_{0z}^2}{p'(V)}\right) + \frac{C(1+N(T))}{L} (\|\phi_0\|^2 + \|\psi_0\|^2) + \frac{C T N^{2}(T)}{L}
				\\
				& + \int_{0}^{t} \int_{-2L}^{2L} \kappa^{2} \left(\frac{(m+1) V_z}{V^{m+2}} r r_{z} - W_{z} r_{z} + \frac{p''(V) V_{z}}{p'(V)} r \psi_z +  G_z \frac{\psi_z}{p'(V)}\right).
			\end{aligned}
		\end{equation*}
		Letting $L \to + \infty$, it holds that
		\begin{equation}
			\begin{aligned}
				& \frac{1}{2} \int \left(r^2 - \frac{\psi_z^2}{p'(V)}\right) + \int_{0}^{t} \int \frac{r_{z}^2}{V^{m+1}} + \int_{0}^{t} \int \frac{s p''(V) V_z}{2 p'(V)^2} \psi_z^2
				\\
				\leq & C (\|r_{0}\|^2 + \|\psi_{0z}\|^2) + \int_{0}^{t} \int \left|\frac{(m+1) V_z}{V^{m+2}} r r_{z} - W_{z} r_{z} + \frac{p''(V) V_{z}}{p'(V)} r \psi_z +  G_z \frac{\psi_z}{p'(V)}\right|.
			\end{aligned}
			\label{3-4411}
		\end{equation}
		By \eqref{3-16}, \eqref{3-18}, $\|r/V^{m+1}\|_{L^\infty}\leq CN(t)$ and $\|r_{z}/V^{m+1}\|_{L^\infty}\leq CN(t)$, Young's inequality gives
		\begin{equation}
			\begin{aligned}
				& \int_{0}^{t} \int \left|\frac{(m+1) V_z}{V^{m+2}} r r_{z} - W_{z} r_{z} + \frac{p''(V) V_{z}}{p'(V)} r \psi_z +  G_z \frac{\psi_z}{p'(V)}\right|
				\\
				\leq & C \int_{0}^{t} \int \left(|r r_{z}| + \frac{r^2 |r_{z}|}{V} + \frac{|r| r_{z}^2}{V^{m+2}} + |V_{z}| |r \psi_z| + |V_z| r^2 |\psi_z| + |r| |r_{z}| |\psi_z|\right)
				\\
				\leq & \left(\frac{1}{2} + C N(t)\right) \int_{0}^{t} \int \frac{r_{z}^2}{V^{m+1}} + C \int_{0}^{t} \int \frac{r^2}{V^{m+1}} + \frac{1}{4} \int_{0}^{t} \int \frac{s p''(V) V_z}{2p'(V)^2} \psi_z^2.
			\end{aligned}
			\label{3-44}
		\end{equation}
		Substituting \eqref{3-22} and \eqref{3-44} into \eqref{3-4411}, we have
		\begin{equation*}
			\begin{aligned}
				\int \left(r^2 - \frac{\psi_z^2}{p'(V)}\right) + \int_0^t \int \frac{r_{z}^2}{V^{m+1}} + \int_0^t \int \frac{s p''(V) V_z}{p'(V)^2} \psi_z^2 \leq C \int (r_{0}^2 + \psi_{0 z}^2 + \phi_0^2 + \psi_0^2).
			\end{aligned}
		\end{equation*}
		We thus finish the proof of \eqref{3-41}.
	\end{proof}
	
	\begin{lemma}\label{lemma3.8}
		Let the assumptions of Theorem \ref{theorem3.1} hold. If $N(T) \ll 1$, then there exists a constant $C>0$ independent of $T$ such that
		\begin{equation}
			\|q_z(t)\|_{w_2}^2+\int_0^t\|q_{zz}\|_{w_3}^2 \mathrm{d}\tau
			\leq C \bigl(\|\phi_0\|^2+\|\psi_0\|_{1}^2
			+\|r_0\|_{w_2}^2+\|q_0\|_{1,w_2}^2\bigr),
			\label{3-53}
		\end{equation}
		for any $t \in[0, T]$.
	\end{lemma}
	\begin{proof}
		We differentiate \eqref{4.7} with respect to $z$ to get 
		\begin{equation}
			\begin{aligned}
				q_{zt} = & sq_{zz} + \left[\frac{q_{z}}{(V+r)^{m+1}} -s\left(\frac{V_z}{(V+r)^{m+1}}-\frac{V_z}{V^{m+1}}\right) - (p(V+r)-p(V))\right]_{zz}
				\\
				= & sq_{zz} + \bigg[ \frac{q_{zz}}{(V+r)^{m+1}} - \frac{(m+1)(V_{z} + r_{z})q_z}{(V+r)^{m+2}} -s\left(\frac{V_{zz}}{(V+r)^{m+1}}-\frac{V_{zz}}{V^{m+1}}\right)
				\\
				& \quad + s (m+1) V_{z}\left(\frac{V_{z} + r_{z}}{(V+r)^{m+2}}-\frac{V_{z}}{V^{m+2}}\right) - (p(V+r)-p(V))_{z}\bigg]_{z}.
			\end{aligned}
			\label{4.8}
		\end{equation}
		Multiplying \eqref{4.8} by $\frac{q_{z}}{V^{2(m+1)}}$, we obtain
		\begin{equation}
			\begin{aligned}
				\frac{1}{2} \frac{\mathrm{d}}{\mathrm{d}t} \frac{q_{z}^{2}}{V^{2(m+1)}} = & \left\{  \left[\frac{q_{zz}}{(V+r)^{m+1}} - \frac{(m+1)(V_{z} + r_{z})q_z}{(V+r)^{m+2}} -s\left(\frac{V_{zz}}{(V+r)^{m+1}}-\frac{V_{zz}}{V^{m+1}}\right) \right. \right.
				\\
				& \left. \left. + s (m+1) V_{z}\left(\frac{V_{z} + r_{z}}{(V+r)^{m+2}}-\frac{V_{z}}{V^{m+2}}\right) - (p(V+r)-p(V))_{z}\right] \frac{q_{z}}{V^{2(m+1)}} \right.
				\\
				& \left. + \frac{s}{2} \frac{q_{z}^{2}}{V^{2(m+1)}} \right\}_{z} + s \frac{(m+1) V_{z} q_{z}^{2}}{V^{2m+3}} - \bigg[ \frac{q_{zz}}{(V+r)^{m+1}} - \frac{(m+1)(V_{z} + r_{z})q_z}{(V+r)^{m+2}} 
				\\
				& -s\left(\frac{V_{zz}}{(V+r)^{m+1}}-\frac{V_{zz}}{V^{m+1}}\right) + s (m+1) V_{z}\left(\frac{V_{z} + r_{z}}{(V+r)^{m+2}}-\frac{V_{z}}{V^{m+2}}\right) 
				\\
				& - (p(V+r)-p(V))_{z}\bigg] \left(\frac{q_{zz}}{V^{2(m+1)}}-\frac{2(m+1)V_{z}q_z}{V^{2m+3}}\right).
			\end{aligned}
			\label{4.9}
		\end{equation}
		Then, multiplying \eqref{4.9} by $\kappa^{2}(z)$ and integrating over $\mathbb{R} \times [0,t]$, we have
		\begin{equation}
			\begin{aligned}
				& \frac{1}{2} \int_{-2L}^{2L} \kappa^{2} \frac{q_{z}^{2}}{V^{2(m+1)}} + \int_{0}^{t} \int_{-2L}^{2L} \frac{\kappa^{2} q_{zz}^{2}}{V^{2(m+1)}(V+r)^{m+1}}
				\\
				= & \frac{1}{2} \int_{-2L}^{2L} \kappa^{2} \frac{q_{0z}^{2}}{V^{2(m+1)}} + s \int_{0}^{t} \int_{-2L}^{2L} \kappa^{2} \frac{(m+1) V_{z} q_{z}^{2}}{V^{2m+3}} + 2 \int_{0}^{t} \int_{-2L}^{2L} \frac{\kappa^{2} (m+1)V_{z} q_{z z} q_{z}}{V^{2m+3}(V+r)^{m+1}}  
				\\
				& + \int_{0}^{t} \int_{-2L}^{2L} \kappa^{2} \left[ \frac{(m+1)(V_{z} + r_{z})q_z}{(V+r)^{m+2}} + s\left(\frac{V_{zz}}{(V+r)^{m+1}}-\frac{V_{zz}}{V^{m+1}}\right) - s (m+1) V_{z} \right.
				\\
				& \left. \quad \times \left(\frac{V_{z} + r_{z}}{(V+r)^{m+2}}-\frac{V_{z}}{V^{m+2}}\right) + (p(V+r)-p(V))_{z}\right] \left(\frac{q_{zz}}{V^{2(m+1)}}-\frac{2(m+1)V_{z}q_z}{V^{2m+3}}\right) 
				\\
				& - \int_{0}^{t} \int_{-2L}^{2L} 2 \kappa \kappa_{z} \left\{  \left[\frac{q_{zz}}{(V+r)^{m+1}} - \frac{(m+1)(V_{z} + r_{z})q_z}{(V+r)^{m+2}} -s\left(\frac{V_{zz}}{(V+r)^{m+1}}-\frac{V_{zz}}{V^{m+1}}\right) \right. \right.
				\\
				& \left. \left. + s (m+1) V_{z}\left(\frac{V_{z} + r_{z}}{(V+r)^{m+2}}-\frac{V_{z}}{V^{m+2}}\right) - (p(V+r)-p(V))_{z}\right] \frac{q_{z}}{V^{2(m+1)}} + \frac{s}{2} \frac{q_{z}^{2}}{V^{2(m+1)}} \right\}.
			\end{aligned}
			\label{4.11}
		\end{equation}
		Now, we estimate the last term of \eqref{4.11}. Applying Young's inequality, Lemmas \ref{lemma3.3}, \ref{lemma3.6} and \ref{lemma3.5}, one gets
		\begin{align*}
			& \bigg| \int_{0}^{t} \int_{-2L}^{2L} 2 \kappa \kappa_{z} \left\{  \left[\frac{q_{zz}}{(V+r)^{m+1}} - \frac{(m+1)(V_{z} + r_{z})q_z}{(V+r)^{m+2}} -s\left(\frac{V_{zz}}{(V+r)^{m+1}}-\frac{V_{zz}}{V^{m+1}}\right) \right. \right.
			\\
			& \left. \left. + s (m+1) V_{z}\left(\frac{V_{z} + r_{z}}{(V+r)^{m+2}}-\frac{V_{z}}{V^{m+2}}\right) - (p(V+r)-p(V))_{z}\right] \frac{q_{z}}{V^{2(m+1)}} + \frac{s}{2} \frac{q_{z}^{2}}{V^{2(m+1)}} \right\} \bigg|
			\\
			\leq & \frac{C}{L} \int_{0}^{t} \int_{-2L}^{2L} \kappa \left(\frac{|q_{z}||q_{zz}|}{V^{3(m+1)}} + \frac{|V_{z}|q_{z}^{2}}{V^{3m+4}} + \frac{|r_{z}| q_{z}^{2}}{V^{3m+4}} + \frac{|V_{zz}||r| |q_{z}|}{V^{3m+4}} + \frac{V_{z}^{2} |r| |q_{z}|}{V^{3m+5}} + \frac{|V_{z}| |r_{z}| |q_{z}|}{V^{3m+4}} \right.
			\\
			& \left. \quad + \frac{|V_{z}| |r| |q_{z}|}{V^{2(m+1)}} + \frac{|r_{z}| |q_{z}|}{V^{2(m+1)}} + \frac{q_{z}^{2}}{V^{2(m+1)}}\right)
			\\
			\leq & \frac{C}{L} \int_{0}^{t} \int_{-2L}^{2L} \left( \kappa^{2} \frac{q_{zz}^{2}}{V^{3(m+1)}} + \frac{q_{z}^{2}}{V^{3(m+1)}} + \frac{q_{z}^{2}}{V^{2(m+1)}}+ \frac{N(t) q_{z}^{2}}{V^{2m+3}}  + \frac{r^{2}}{V^{m+1}} + \frac{r_{z}^{2}}{V^{m+1}} \right)
			\\
			\leq & \frac{C}{L} \int_{0}^{t} \int_{-2L}^{2L} \kappa^{2} \frac{q_{zz}^{2}}{V^{3(m+1)}} + \frac{C}{L} (\|\phi_0\|^2 + \|\psi_{0}\|^2 + \|r_{0}\|_{w_2}^2 + \|q_{0}\|_{w_{2}}^2),
		\end{align*}
		where we have used the estimates \eqref{3-16}, $\|r_{z}/V^{m+1}\|_{L^{\infty}} \leq C N(T)$, $c_{1}/V^{m+1} \leq 1/(V+r)^{m+1} \leq C/V^{m+1}$, $1 \leq C/V^{m}$ and \eqref{4.10}. Based on the above inequality, \eqref{4.11} becomes
		\begin{align*}
			& \frac{1}{2} \int_{-2L}^{2L} \kappa^{2} \frac{q_{z}^{2}}{V^{2(m+1)}} + \left(c_{1} - \frac{C}{L}\right) \int_{0}^{t} \int_{-2L}^{2L} \frac{\kappa^{2} q_{zz}^{2}}{V^{3(m+1)}}
			\\
			\leq & \frac{1}{2} \int_{-2L}^{2L} \kappa^{2} \frac{q_{0z}^{2}}{V^{2(m+1)}} + \frac{C}{L} (\|\phi_0\|^2 + \|\psi_{0}\|^2 + \|r_{0}\|_{w_2}^2 + \|q_{0}\|_{w_{2}}^2) + s \int_{0}^{t} \int_{-2L}^{2L} \kappa^{2} \frac{(m+1) V_{z} q_{z}^{2}}{V^{2m+3}}   
			\\
			& + 2 \int_{0}^{t} \int_{-2L}^{2L} \frac{\kappa^{2} (m+1)V_{z} q_{z z} q_{z}}{V^{2m+3}(V+r)^{m+1}} + \int_{0}^{t} \int_{-2L}^{2L} \kappa^{2} \left[ \frac{(m+1)(V_{z} + r_{z})q_z}{(V+r)^{m+2}} + s\left(\frac{V_{zz}}{(V+r)^{m+1}}-\frac{V_{zz}}{V^{m+1}}\right) \right.
			\\
			& \left. \quad - s (m+1) V_{z} \left(\frac{V_{z} + r_{z}}{(V+r)^{m+2}}-\frac{V_{z}}{V^{m+2}}\right) + (p(V+r)-p(V))_{z}\right] \left(\frac{q_{zz}}{V^{2(m+1)}}-\frac{2(m+1)V_{z}q_z}{V^{2m+3}}\right).
		\end{align*}
		Letting $L \to +\infty$, we derive from $1/(V+r)^{m+1} \leq C/V^{m+1}$ that 
		\begin{equation}
			\begin{aligned}
				& \frac{1}{2} \int \frac{q_{z}^{2}}{V^{2(m+1)}} + \frac{c_{1}}{2} \int_{0}^{t} \int \frac{q_{zz}^{2}}{V^{3(m+1)}}
				\\
				\leq & C \int \frac{q_{0z}^{2}}{V^{2(m+1)}} + C \int_{0}^{t} \int \frac{ |V_{z}| q_{z}^{2}}{V^{2m+3}} + C \int_{0}^{t} \int \frac{ |V_{z}| |q_{z z}| |q_{z}|}{V^{3m+4}} + C \int_{0}^{t} \int \bigg| \frac{(m+1)(V_{z} + r_{z})q_z}{(V+r)^{m+2}} 
				\\
				& \quad + s\left(\frac{V_{zz}}{(V+r)^{m+1}}-\frac{V_{zz}}{V^{m+1}}\right) - s (m+1) V_{z} \left(\frac{V_{z} + r_{z}}{(V+r)^{m+2}}-\frac{V_{z}}{V^{m+2}}\right) 
				\\
				& \quad + (p(V+r)-p(V))_{z}\bigg| \left(\frac{|q_{zz}|}{V^{2(m+1)}}+\frac{2(m+1)|V_{z}||q_z|}{V^{2m+3}}\right).
			\end{aligned}
			\label{4.13}
		\end{equation}
		Thanks to $1 \leq C/V$, \eqref{3-16}, \eqref{rz:infinity}, we use Taylor's expansion, Lemmas \ref{lemma3.3}, \ref{lemma3.6} and \ref{lemma3.5} to obtain
		\begin{align*}
			\int_{0}^{t} \int \frac{ |V_{z}| q_{z}^{2}}{V^{2m+3}} + \int_{0}^{t} \int \frac{ |V_{z}| |q_{z z}| |q_{z}|}{V^{3m+4}} \leq & C \int_{0}^{t} \int \frac{q_{z}^{2}}{V^{m+1}} + \frac{c_{1}}{8} \int_{0}^{t} \int \frac{q_{zz}^{2}}{V^{3(m+1)}} 
			\\
			\leq & C (\|\phi_0\|^2 + \|\psi_0\|^2 + \|r_{0}\|_{w_2}^2 + \|q_{0}\|_{w_{2}}^2) + \frac{c_{1}}{8} \int_{0}^{t} \int \frac{q_{zz}^{2}}{V^{3(m+1)}},
		\end{align*}
		and 
		\begin{align*}
			& \int_{0}^{t} \int \bigg| \frac{(m+1)(V_{z} + r_{z})q_z}{(V+r)^{m+2}} + s\left(\frac{V_{zz}}{(V+r)^{m+1}}-\frac{V_{zz}}{V^{m+1}}\right) - s (m+1) V_{z} \left(\frac{V_{z} + r_{z}}{(V+r)^{m+2}}-\frac{V_{z}}{V^{m+2}}\right) 
			\\
			& \quad + (p(V+r)-p(V))_{z}\bigg| \left(\frac{|q_{zz}|}{V^{2(m+1)}}+\frac{2(m+1)|V_{z}||q_z|}{V^{2m+3}}\right)
			\\
			\leq & C \int_{0}^{t} \int \bigg[\frac{|V_{z}||q_z|}{V^{m+2}} + \frac{|r_{z}||q_z|}{V^{m+2}} + \frac{|V_{zz}||r|}{V^{m+2}} + \frac{V_{z}^{2}|r|}{V^{m+3}} + \frac{|V_{z}||r_z|}{V^{m+2}} + |r_{z}| + |V_{z}| |r|\bigg] \left(\frac{|q_{zz}|}{V^{2(m+1)}}+\frac{|V_{z}||q_z|}{V^{2m+3}}\right)
			\\
			\leq & C \int_{0}^{t} \int \left(\frac{q_{z}^{2}}{V^{m+1}} + \frac{r^{2}}{V^{m+1}} + \frac{r_{z}^{2}}{V^{m+1}}\right) + \left(\frac{c_{1}}{8} + CN(t)\right) \int_{0}^{t} \int \frac{q_{zz}^{2}}{V^{3(m+1)}} + CN(t)\int_{0}^{t} \int \frac{q_{z}^{2}}{V^{3(m+1)}} 
			\\
			\leq & C (\|\phi_0\|^2 + \|\psi_0\|_{1}^2 + \|r_{0}\|_{w_2}^2 + \|q_{0}\|_{w_{2}}^2) + \left(\frac{c_{1}}{8} + CN(t)\right) \int_{0}^{t} \int \frac{q_{zz}^{2}}{V^{3}}.
		\end{align*}
		Combining the above estimates, \eqref{4.13} takes the form
		\begin{equation*}
			\frac{1}{2} \int \frac{q_{z}^{2}}{V^{2(m+1)}} + \frac{c_{1}}{8} \int_{0}^{t} \int \frac{q_{zz}^{2}}{V^{3(m+1)}} \leq C (\|\phi_0\|^2 + \|\psi_0\|_{1}^2 + \|r_{0}\|_{w_2}^2 + \|q_{0}\|_{1, w_{2}}^2),
		\end{equation*}
		provided $N(t)$ is suitably small. 
	\end{proof}
	
	\begin{lemma}\label{lemma3.9}
		Let the assumptions of Theorem \ref{theorem3.1} hold. If $N(T) \ll 1$, then there exists a constant $C>0$ independent of $T$ such that
		\begin{equation}
			\|r_z(t)\|_{w_2}^2 \leq C (\|\phi_0\|^2 + \|\psi_0\|_{1}^2 + \|r_{0}\|_{1,w_2}^2 + \|q_{0}\|_{1, w_{2}}^2),
			\label{3-57}
		\end{equation}
		for any $t \in[0, T]$.
	\end{lemma}
	\begin{proof}
		Differentiating $r_t=s r_z+q_z$ with respect to $z$, multiplying it by $\kappa^{2} r_z/V^{2(m+1)}$ and integrating over $\mathbb{R} \times [0,t]$, we obtain
		\begin{equation*}
			\begin{aligned}
				\frac{1}{2} \int_{-2L}^{2L} \frac{\kappa^{2} r_z^2}{V^{2(m+1)}} = & \frac{1}{2} \int_{-2L}^{2L} \frac{\kappa^{2} r_{0z}^{2}}{V^{2(m+1)}} + \int_{0}^{t} \int_{-2L}^{2L} \frac{\kappa^{2} q_{zz}r_z}{V^{2(m+1)}} 
				\\
				& + s \int_{0}^{t} \int_{-2L}^{2L} \kappa^{2} \frac{(m+1) V_z r_z^2}{V^{2m+3}} - s \int_{0}^{t} \int_{-2L}^{2L} \frac{\kappa \kappa_{z} r_{z}^{2}}{V^{2(m+1)}}.
			\end{aligned}
		\end{equation*}
		By \eqref{3-11}, we have
		\begin{equation*}
			\left| - s \int_{0}^{t} \int_{-2L}^{2L} \frac{\kappa \kappa_{z} r_{z}^{2}}{V^{2(m+1)}} \right| \leq \frac{CTN^{2}(T)}{L}.
		\end{equation*}
		Passing to the limit $L \to + \infty$ and using Young's inequality, we get
		\begin{equation*}
			\begin{aligned}
				\frac{1}{2} \int \frac{r_z^2}{V^{2(m+1)}} \leq & \frac{1}{2} \int \frac{r_{0z}^{2}}{V^{2(m+1)}} + \left|\int_{0}^{t} \int \frac{q_{zz}r_z}{V^{2(m+1)}}\right| + \left|s \int_{0}^{t} \int \frac{(m+1)V_zr_z^2}{V^{2m+3}}\right|
				\\
				\leq & \frac{1}{2} \int \frac{r_{0z}^{2}}{V^{2(m+1)}} + C \int_{0}^{t} \int \frac{q_{zz}^{2}}{V^{3(m+1)}} + C \int_{0}^{t} \int \frac{r_{z}^{2}}{V^{m+1}}.
			\end{aligned}
		\end{equation*}
		Applying Lemmas \ref{lemma3.5} and \ref{lemma3.8}, we complete the proof.
	\end{proof}
	
	\begin{lemma}\label{lemma3.11}
		Let the assumptions of Theorem \ref{theorem3.1} hold. If $N(T) \ll 1$, then there exists a constant $C>0$ independent of $T$ such that
		\begin{equation}
			\|r_{z}\|^2 + \|\psi_{zz}\|^2 + \int_0^t \|r_{zz}\|_{w_1}^2  + \int_0^t \int d(V) \psi_{zz}^2 \leq C (\|\phi_0\|_2^2 + \|\psi_0\|_2^2),
			\label{3-91}
		\end{equation}
		for any $t \in [0, T]$.
	\end{lemma}
	\begin{proof}
		Differentiating \eqref{3-42} with respect to $z$, we get
		\begin{equation}
			\begin{cases}
				r_{zt} = sr_{zz} + \psi_{zzz} + (r/V^{m+1})_{zzz} + W_{zzz},
				\\
				\psi_{zzt} = s\psi_{zzz} - p'(V) r_{zz} - 2p''(V)V_zr_z - p'''(V)V_z^2 r - p''(V)V_{zz}r - G_{zz}.
			\end{cases}
			\label{3-421}
		\end{equation}
		Multiplying $\eqref{3-421}_{1}$ by $r_{z}$ and $\eqref{3-421}_{2}$ by $(-\frac{\psi_{zz}}{p'(V)})$, respectively, and then adding them gives
		\begin{align*}
			\frac{1}{2} \frac{\mathrm{d}}{\mathrm{d}t} \left(r_{z}^2 - \frac{\psi_{zz}^2}{p'(V)}\right) = & \left[ \left(\frac{r_{zz}}{V^{m+1}} - \frac{2 (m+1) V_{z} r_{z}}{V^{m+2}} - \frac{(m+1) V_{zz} r}{V^{m+2}} + \frac{(m+1)(m+2)V_{z}^{2} r}{V^{m+3}}\right) r_{z} \right.
			\\
			& \left. + \frac{s}{2} r_{z}^{2} + \psi_{zz} r_{z} + W_{zz} r_{z} + \frac{s\psi_{zz}^2}{-2p'(V)}\right]_{z} -  \frac{r_{zz}^2}{V^{m+1}} + \frac{2 (m+1) V_{z} }{V^{m+2}} r_{z} r_{zz}  
			\\
			& + \frac{(m+1) V_{zz}}{V^{m+2}} r r_{zz} - \frac{(m+1)(m+2)V_{z}^{2}}{V^{m+3}} r r_{zz}- \frac{s p''(V) V_z}{2 p'(V)^2} \psi_{zz}^2 - W_{zz} r_{zz}   
			\\
			& + \frac{2p''(V)V_z r_z \psi_{zz}}{p'(V)} + \frac{ p'''(V)V_z^2 r \psi_{zz} }{p'(V)} + \frac{ p''(V)V_{zz} r \psi_{zz}}{p'(V)} +  G_{zz} \frac{\psi_{zz}}{p'(V)}.
		\end{align*}
		Multiplying it by $\kappa^{2}(z)$ and integrating over $\mathbb{R} \times [0,t]$, one can show that
		\begin{equation}
			\begin{aligned}
				& \frac{1}{2} \int_{-2L}^{2L} \kappa^{2} \left(r_{z}^2 - \frac{\psi_{zz}^2}{p'(V)}\right) + \int_{0}^{t} \int_{-2L}^{2L} \kappa^{2} \frac{r_{zz}^2}{V^{m+1}} + \int_{0}^{t} \int_{-2L}^{2L} \kappa^{2} \frac{s p''(V) V_z}{2 p'(V)^2} \psi_{zz}^2
				\\
				= & \frac{1}{2} \int_{-2L}^{2L} \kappa^{2} \left(r_{0z}^2 - \frac{\psi_{0zz}^2}{p'(V)}\right) + \int_{0}^{t} \int_{-2L}^{2L} \kappa^{2} \bigg[\frac{2 (m+1) V_{z} }{V^{m+2}} r_{z} r_{zz} + \frac{(m+1) V_{zz}}{V^{m+2}} r r_{zz}  
				\\
				& - \frac{(m+1)(m+2)V_{z}^{2}}{V^{m+3}} r r_{zz} - W_{zz} r_{zz} + \frac{2p''(V)V_z r_z \psi_{zz}}{p'(V)} + \frac{ p'''(V)V_z^2 r \psi_{zz} }{p'(V)} 
				\\
				& + \frac{ p''(V)V_{zz} r \psi_{zz}}{p'(V)} +  G_{zz} \frac{\psi_{zz}}{p'(V)}\bigg] - \int_{0}^{t} \int_{-2L}^{2L} 2 \kappa \kappa_{z} \bigg[ \left(\frac{r_{zz}}{V^{m+1}} - \frac{2 (m+1) V_{z} r_{z}}{V^{m+2}} \right.
				\\
				& \left. \quad - \frac{(m+1) V_{zz} r}{V^{m+2}} + \frac{(m+1)(m+2)V_{z}^{2} r}{V^{m+3}}\right) r_{z} + \frac{s}{2} r_{z}^{2} + \psi_{zz} r_{z} + W_{zz} r_{z} + \frac{s\psi_{zz}^2}{-2p'(V)}\bigg].
			\end{aligned}
			\label{3-451}
		\end{equation}
		We utilize the estimates $\|\psi_{zz}\| \leq N(t)$, $\|r/V\|_{L^{\infty}},\|r_{z}/V\|_{L^{\infty}} \leq C N(t)$, together with \eqref{3-16}, \eqref{3-18}, Lemmas \ref{lemma3.3} and \ref{lemma3.5}, to get 
		\begin{align*}
			& \left|\int_{0}^{t} \int_{-2L}^{2L} 2 \kappa \kappa_{z} \bigg[\left(\frac{r_{zz}}{V^{m+1}} - \frac{2 (m+1) V_{z} r_{z}}{V^{m+2}} - \frac{(m+1) V_{zz} r}{V^{m+2}} + \frac{(m+1)(m+2)V_{z}^{2} r}{V^{m+3}}\right) r_{z} \right.
			\\
			& \left. \quad + \frac{s}{2} r_{z}^{2} + \psi_{zz} r_{z} + W_{zz} r_{z} + \frac{s\psi_{zz}^2}{-2p'(V)}\bigg]\right|
			\\
			\leq & \frac{C}{L} \int_{0}^{t} \int_{-2L}^{2L} \kappa \bigg[\frac{|r_{zz}| |r_{z}|}{V^{m+1}} + r_{z}^{2} + V^{m+1}|r||r_{z}| + |\psi_{zz}| |r_{z}| + V^{m} |r_{z}| |r|^{2}
			\\
			& \quad + \frac{|r_{z}|^{3}}{V^{m+2}} + \frac{|r_{z}|^{2} |r|}{V} + \frac{|r||r_{z}||r_{zz}|}{V^{m+2}} + \psi_{zz}^{2}\bigg]
			\\
			\leq & \frac{C (1+N(T))}{L} \int_{0}^{t} \int_{-2L}^{2L} \kappa^{2} \frac{r_{zz}^{2}}{V^{m+1}} + \frac{C (1+N(T))}{L} \int_{0}^{t} \int_{-2L}^{2L} \left(\frac{r^{2}}{V^{m+1}} + \frac{r_{z}^{2}}{V^{m+1}}\right) + \frac{C}{L} \int_{0}^{t} \int_{-2L}^{2L} \psi_{zz}^{2}
			\\
			\leq & \frac{C (1+N(T))}{L} \int_{0}^{t} \int_{-2L}^{2L} \kappa^{2} \frac{r_{zz}^{2}}{V^{m+1}} + \frac{C (1+N(T))}{L} (\|\phi_0\|_{1}^2 + \|\psi_0\|_{1}^2) + \frac{C T N^{2}(T)}{L}.
		\end{align*}
		Therefore, \eqref{3-451} can be estimated as
		\begin{align*}
			& \frac{1}{2} \int_{-2L}^{2L} \kappa^{2} \left(r_{z}^2 - \frac{\psi_{zz}^2}{p'(V)}\right) + \left(1 - \frac{C (1+N(T))}{L}\right) \int_{0}^{t} \int_{-2L}^{2L} \kappa^{2} \frac{r_{zz}^2}{V^{m+1}} + \int_{0}^{t} \int_{-2L}^{2L} \kappa^{2} \frac{s p''(V) V_z}{2 p'(V)^2} \psi_{zz}^2
			\\
			\leq & \frac{1}{2} \int_{-2L}^{2L} \kappa^{2} \left(r_{0z}^2 - \frac{\psi_{0zz}^2}{p'(V)}\right) + \frac{C (1+N(T))}{L} (\|\phi_0\|_{1}^2 + \|\psi_0\|_{1}^2) + \frac{C T N^{2}(T)}{L}
			\\
			& + \int_{0}^{t} \int_{-2L}^{2L} \kappa^{2} \left[\frac{2 (m+1) V_{z} }{V^{m+2}} r_{z} r_{zz} + \frac{(m+1) V_{zz}}{V^{m+2}} r r_{zz} - \frac{(m+1)(m+2)V_{z}^{2}}{V^{m+3}} r r_{zz} - W_{zz} r_{zz} \right.
			\\
			& \left. \quad + \frac{2p''(V)V_z r_z \psi_{zz}}{p'(V)} + \frac{ p'''(V)V_z^2 r \psi_{zz} }{p'(V)} + \frac{ p''(V)V_{zz} r \psi_{zz}}{p'(V)} +  G_{zz} \frac{\psi_{zz}}{p'(V)}\right].
		\end{align*}
		Letting $L \to + \infty$, we deduce that
		\begin{equation*}
			\begin{aligned}
				& \frac{1}{2} \int \left(r_{z}^2 - \frac{\psi_{zz}^2}{p'(V)}\right) + \int_{0}^{t} \int \frac{r_{zz}^2}{V^{m+1}} + \int_{0}^{t} \int \frac{s p''(V) V_z}{2 p'(V)^2} \psi_{zz}^2
				\\
				\leq & C (\|r_{0z}\|^2 + \|\psi_{0zz}\|^2) + \int_{0}^{t} \int \left|\frac{2 (m+1) V_{z} }{V^{m+2}} r_{z} r_{zz} + \frac{(m+1) V_{zz}}{V^{m+2}} r r_{zz} - \frac{(m+1)(m+2)V_{z}^{2}}{V^{m+3}} r r_{zz} \right.
				\\
				& \left. \quad  - W_{zz} r_{zz} + \frac{2p''(V)V_z r_z \psi_{zz}}{p'(V)} + \frac{ p'''(V)V_z^2 r \psi_{zz} }{p'(V)} + \frac{ p''(V)V_{zz} r \psi_{zz}}{p'(V)}\right| +  \int_{0}^{t} \int \left|G_{zz} \frac{\psi_{zz}}{p'(V)}\right|.
			\end{aligned}
		\end{equation*}
		Employing $|V_{z z}| \leq CV|V_z|$, $\|r/V\|_{L^{\infty}} \leq C N(t)$, $\|r_{z}/V\|_{L^{\infty}} \leq C N(t)$, $\|\psi_{zz}\| \leq C N(t)$ and Young's inequality, we have
		\begin{align*}
			& \int_{0}^{t} \int \left|\frac{2 (m+1) V_{z} }{V^{m+2}} r_{z} r_{zz} + \frac{(m+1) V_{zz}}{V^{m+2}} r r_{zz} - \frac{(m+1)(m+2)V_{z}^{2}}{V^{m+3}} r r_{zz} - W_{zz} r_{zz}  \right.
			\\
			& \left. \quad  + \frac{2p''(V)V_z r_z \psi_{zz}}{p'(V)} + \frac{ p'''(V)V_z^2 r \psi_{zz} }{p'(V)} + \frac{ p''(V)V_{zz} r \psi_{zz}}{p'(V)} \right|
			\\
			\leq & C \int_{0}^{t} \int \left(|r_{z} r_{zz}| + V^{m+1}|r r_{zz}| + V^{m} r^{2}| r_{zz}| + \frac{|r_{zz}| r_{z}^2}{V^{m+2}} + \frac{|r| |r_{z}||r_{zz}| }{V} \right. 
			\\
			& \left. \quad + \frac{|r| r_{zz}^{2}}{V^{m+2}} + |V_{z}| |r_{z} \psi_{zz}| + V |V_{z}| |r \psi_{zz}|\right)
			\\
			\leq & \left(\frac{1}{8} + C N(t)\right) \int_{0}^{t} \int \frac{r_{zz}^2}{V^{m+1}} + C (1+N(t)) \int_{0}^{t} \int \left(\frac{r^2}{V^{m+1}} + \frac{r_{z}^2}{V^{m+1}}\right) + \frac{1}{4} \int_{0}^{t} \int \frac{s p''(V) V_z}{2p'(V)^2} \psi_{zz}^2,
		\end{align*}
		and
		\begin{align*}
			\int_{0}^{t} \int \left| G_{zz} \frac{\psi_{zz}}{p'(V)} \right| \leq & C \int_{0}^{t} \int \left(V_z^2r^2 |\psi_{zz}| + |V_z||r||r_{z}||\psi_{zz}| + |V_{zz}|r^2 |\psi_{zz}| + r_{z}^2 |\psi_{zz}| + |r||r_{zz}| |\psi_{zz}|\right)
			\\
			\leq & C \int_{0}^{t} \left( \|r\|_{L^{\infty}} \|r\| + \|r\|_{L^{\infty}} \|r_{z}\| + \|r_{z}\|_{L^{\infty}} \|r_{z}\| + \|r\|_{L^{\infty}} \|r_{zz}\|\right) \|\psi_{zz}\|
			\\
			\leq & C N(t) \int_{0}^{t} \left( \|r_{z}\|^{\frac{1}{2}} \|r\|^{\frac{3}{2}} + \|r_{z}\|^{\frac{3}{2}} \|r\|^{\frac{1}{2}} + \|r_{zz}\|^{\frac{1}{2}} \|r_{z}\|^{\frac{3}{2}} + \|r\|^{\frac{1}{2}} \|r_{z}\|^{\frac{1}{2}} \|r_{zz}\|\right) 
			\\
			\leq & C N(t) \int_{0}^{t} \left( \|r\|^{2} + \|r_{z}\|^{2} + \|r_{zz}\|^{2}\right).
		\end{align*}
		Using Lemmas \ref{lemma3.3} and \ref{lemma3.5}, it follows from the estimates above that \eqref{3-91}. we complete the proof of Lemma \ref{lemma3.11}.
	\end{proof}
	
	\begin{proof}[Proof of Proposition \ref{proposition3.2}]
		Combining Lemmas \ref{lemma3.3}--\ref{lemma3.11}, we choose $\delta_1 \leq \delta_0$ sufficiently small. If $N(t) \leq \delta_1$ for any $0 \leq t \leq T$, then estimate \eqref{3-12} follows immediately. Hence Proposition \ref{proposition3.2} is proved.
	\end{proof}

	\subsection{Global existence and convergence}
	\begin{proof}[Proof of Theorem \ref{theorem3.1}]
		Let $\delta_2 = \min\{\delta_1/(2C_{*}), \delta_1/(C C_{*})\}$ be sufficiently small. By the standard continuation arguments associated with Propositions \ref{proposition3.1} and \ref{proposition3.2}, we obtain a unique global solution $(\phi, \psi) \in X(0, \infty)$ satisfying estimate \eqref{3-12} for all $t \in[0, \infty)$.
		
		Before proving \eqref{3-13}, we recall a classical result for later use; see \cite{MANKG}.
		\begin{lemma}\label{lemma3.10}
			If $F\in W^{1,1}(0,\infty)$, then $F(t)\to0$ as $t\to\infty$.
		\end{lemma}
		\begin{proof}
			Since $F'\in L^1$, $F$ has a finite limit at infinity.
			Since $F\in L^1$, the limit must be zero.
		\end{proof}
		
		To complete the proof of Theorem \ref{theorem3.1}, utilizing the established estimate \eqref{3-12} for all $t \in[0, \infty)$, we have
		\begin{equation}
			\|r\|_{w_{2}}^2 + \|r_{z}\|_{w_{2}}^2 + \|q\|_{w_{2}}^2 + \|q_{z}\|_{w_{2}}^2 + \int_0^{\infty} (\|r\|_{w_{1}}^2 + \|r_{z}\|_{w_{1}}^2 + \|q\|_{w_{1}}^2 + \|q_{z}\|_{1,w_{3}}^2) \mathrm{d} t \leq C.
			\label{3-612}
		\end{equation}
		In particular, $\|r(\cdot, t)\|^2 \in L^1(0, \infty)$. We next show that $\frac{\mathrm{d}}{\mathrm{d}t}\|r(\cdot, t)\|^2 \in L^1(0, \infty)$. We apply the Cauchy--Schwarz inequality and the identity $r_{t} = sr_{z} + q_{z}$ of \eqref{3-251} to get
		\begin{align*}
			\int_0^{\infty}\left|\frac{\mathrm{d}}{\mathrm{d}t} \|r(\cdot, t)\|^2\right| \mathrm{d}t = \int_0^{\infty}\left|\frac{\mathrm{d}}{\mathrm{d}t} \int r^2 \right| \mathrm{d}t = 2 \int_0^{\infty}\left|\int r r_{t} \right| \mathrm{d}t \leq C \int_0^{\infty} \int (r^2 + r_{z}^2 + q_{z}^2) \mathrm{d}t \leq C.
		\end{align*}
		It follows from Lemma \ref{lemma3.10} that
		\begin{align*}
			\|r(\cdot, t)\|^2 \to 0, \quad \text { as } t \to +\infty. 
		\end{align*}
		This along with \eqref{3-612}, for all $z \in \mathbb{R}$, implies that
		\begin{align*}
			r^2(z, t) = 2 \int_{-\infty}^z r r_{z}(y, t) \mathrm{d}y \leq C \left(\int r^2 \mathrm{d}y\right)^{1 / 2}\left(\int r_{z}^2 \mathrm{d}y\right)^{1 / 2} \leq C \left\|r(\cdot, t)\right\| \to 0, \, \text { as } t \to +\infty,
		\end{align*}
		namely,
		$$
		\sup _{z \in \mathbb{R}}|r(z, t)| \to 0, \quad \text { as } t \to + \infty.
		$$
		Recall from \eqref{4.7} that
		\begin{equation*}
			q_t = sq_z  - (p(V+r)-p(V))_{z} + \left[\frac{q_{z}}{(V+r)^{m+1}} -s\left(\frac{V_z}{(V+r)^{m+1}}-\frac{V_z}{V^{m+1}}\right) \right]_z.
		\end{equation*}
		By \eqref{3-612} and Taylor's expansion, one has
		\begin{align*}
			\int_{0}^{\infty} \int (s q_z)^2 \leq C,
		\end{align*}
		and
		\begin{align*}
			\int_0^\infty \int [-(p(V+r)-p(V))_z]^2 = & \int_0^\infty \int [p'(V+r)r_z + V_z (p'(V+r)-p'(V))]^2
			\\
			\leq & C \int_0^\infty \int (r_{z}^{2} + V_z^2 r^2) 
			\\
			\leq & C \int_0^\infty \left(\int \frac{r_{z}^2}{V^{m+1}}+\int \frac{r^2}{V^{m+1}}\right)
			\\
			\leq & C (\|\phi_0\|_1^2 + \|\psi_0\|_1^2).
		\end{align*}
		Moreover, we use \eqref{3-16}, \eqref{4.711}, \eqref{3-612}, $\| r_{z}/V^{m+1} \|_{L^{\infty}} \leq C N(t)$ and Taylor's expansion to derive
		\begin{align*}
			& \int_0^\infty \int \left[\left(\frac{q_{z}}{(V+r)^{m+1}}\right)_z\right]^{2} + \int_0^{\infty} \int \left[s\left(\frac{V_z}{(V+r)^{m+1}}-\frac{V_z}{V^{m+1}}\right)_z\right]^2
			\\
			= & \int_0^\infty \int \left(\frac{q_{zz}}{(V+r)^{m+1}} - \frac{(m+1)q_{z}(V_{z} + r_{z})}{(V + r)^{m+2}}\right)^{2} + \int_0^\infty \int \left[ s \left(\frac{V_{zz}}{(V+r)^{m+1}}-\frac{V_{zz}}{V^{m+1}}\right) \right.
			\\
			& \left. \quad - s (m+1) V_{z} \left(\frac{V_{z} + r_{z}}{(V+r)^{m+2}}-\frac{V_{z}}{V^{m+2}}\right) \right]^{2}
			\\
			\leq & C \int_0^\infty \int \left(\frac{q_{zz}^{2}}{(V+r)^{2(m+1)}} + \frac{|q_{z} q_{zz} (V_{z} + r_{z})|}{(V+r)^{2m+3}} + \frac{q_{z}^{2} (V_{z} + r_{z})^{2}}{(V+r)^{2m+4}}\right)
			\\
			& + C \int_0^\infty \int \left(\frac{\left|V_{z z}\right|}{V^{m+2}}\left|r\right|+\frac{V_z^2}{V^{m+3}}\left|r\right|+\frac{\left|V_z\right|}{V^{m+2}}\left|r_{z}\right|\right)^2
			\\
			\leq & C \int_0^\infty \int \left(\frac{q_{zz}^{2}}{V^{2(m+1)}} + \frac{|q_{z}q_{zz}|}{V^{m+1}} + \frac{|q_{z}q_{zz}r_{z}|}{V^{2m+3}} + q_{z}^{2} + \frac{q_{z}^{2} r_{z}^{2}}{V^{2m+4}} + V^{2(m+1)} r^{2} + r_{z}^{2}\right)
			\\
			\leq & C \int_0^\infty \int \left(\frac{q_{zz}^{2}}{V^{2(m+1)}} + \frac{q_{z}^{2}}{V^{2(m+1)}} + r_{z}^{2} + r^{2}\right)
			\\
			\leq & C.
		\end{align*}
		Applying the above estimates, we obtain
		\begin{align*}
			\int_0^{\infty} \int q_{t}^{2} \leq C.
		\end{align*}
		Therefore, combining this estimate with \eqref{3-612}, a direct calculation yields
		\begin{align*}
			\int_0^{\infty} \left|\frac{\mathrm{d}}{\mathrm{d}t} \|q(\cdot, t)\|^2\right| \mathrm{d}t  = \int_0^{\infty} \left|\frac{\mathrm{d}}{\mathrm{d}t} \int q^2 \right| \mathrm{d}t = 2 \int_0^{\infty} \left|\int q q_{t}\right| \mathrm{d}t \leq C \int_0^{\infty} \int (q^2 + q_{t}^2) \mathrm{d}t \leq C.
		\end{align*}
		We get from Lemma \ref{lemma3.10} that
		\begin{align}
			\|q(\cdot, t)\|^2 \to 0, \, \text { as } t \to +\infty.
			\label{4-6}
		\end{align}
		Therefore, for all $z \in \mathbb{R}$, by \eqref{3-612}, \eqref{4-6} and the Sobolev inequality, we deduce that
		\begin{align*}
			q^2(z, t) = 2 \int_{-\infty}^z q q_{z}(y, t) \mathrm{d}y \leq C \left(\int q^2 \mathrm{d}y\right)^{1 / 2} \left(\int q_{z}^2 \mathrm{d}y\right)^{1 / 2} \leq C \|q(\cdot, t)\| \to 0, \, \text { as } t \to +\infty,
		\end{align*}
		namely,
		$$
		\sup _{z \in \mathbb{R}}|q(z, t)| \to 0, \quad \text { as } t \to + \infty.
		$$
		Now we complete the proof of Theorem \ref{theorem3.1}.
	\end{proof}

	\begin{proof}[Proof of Theorem \ref{theorem1.3}]
		First take $x_0=0$. it follows from the definitions of $b_0$ and $B$ that
		\begin{equation*}
			q_0=\psi_{0z} +\frac{v_{0z}}{v_0^{m+1}} -\frac{V_z}{V^{m+1}} =u_0-U.
		\end{equation*}
		Thus the smallness assumption implies the hypotheses of
		Theorem \ref{theorem3.1}. We reconstruct
		\begin{equation*}
			v(z+st,t)=V(z)+r(z,t),\qquad
			u(z+st,t)=U(z)+q(z,t).
		\end{equation*}
		By $r_t=s r_z+q_z$ and $U_z=-sV_z$, we obtain the mass equation. \eqref{4.7} gives
		\begin{align*}
			q_t = sq_z  - (p(V+r)-p(V))_{z} + \left(\frac{u_{z} }{v^{m+1}} - \frac{U_z}{V^{m+1}}\right)_z.
		\end{align*}
		From \eqref{1-111}, \eqref{1-5}, \eqref{1-10}, \eqref{2-15} and the preceding identity, we get the momentum equation
		\begin{align*}
			u_t-s u_z+p(v)_z = \left(\frac{u_z}{v^{m+1}}\right)_z.
		\end{align*}
		The regularity and uniqueness of the solution follow from the perturbation theorem and $V+r \geq V/2>0$. The convergence is an immediate consequence of \eqref{3-13}. Finally, restoring the shift $x_0$ completes the proof of the theorem.
	\end{proof}

	\noindent {\bf Declaration}. This is an original study, no AI is applied for mathematical ideas and proofs.

\end{document}